\documentclass[a4paper,10pt]{amsart}

\usepackage[latin1]{inputenc}
\usepackage{amssymb,amsmath}

   \usepackage[pagebackref]{hyperref}
  \hypersetup{
colorlinks=true,
linkcolor=cyan,
citecolor=cyan
}
 
\usepackage{mathtools}

\usepackage{color}

\theoremstyle{definition}

\numberwithin{equation}{section}

\def\elem(#1,#2){  \{ \frac{#1}  {\overline {\ #2\ }} \} }

\title[From singularities to graphs]{From singularities to graphs}
\author{Patrick Popescu-Pampu}
   \address{Univ. Lille, CNRS, UMR 8524 - Laboratoire Paul Painlev\'e, F-59000 Lille, 
         France.}
   \email{patrick.popescu-pampu@univ-lille.fr}
\date{21 June 2022}

\subjclass[2010]{14B05 (primary), 32S25, 32S45, 32S50, 57M15, 01A60}
\keywords{Coxeter diagrams, Dual graphs, Du Val singularities, Graph manifolds, 
     Models of surfaces, Plumbing calculus, 
    Resolution of singularities, Singularity links, Surface singularities.}

\begin{document}
%\linenumbers

\begin{abstract}
      In this text I present some problems which led to the introduction of special kinds of 
      graphs as tools for studying singular points of algebraic surfaces. I explain how  
      such graphs were first described using words, and how  
      several classification problems made it necessary to draw them, leading  
      to the elaboration of a special kind of calculus with graphs. 
      This non-technical paper is intended to be readable both by mathematicians and 
      philosophers or historians of mathematics. 
\end{abstract}

{\bf With some differences in the management of references and section numbering, this paper will appear as a chapter of the book {\em When Form Becomes Substance. Power of Gestures, Diagrammatical Intuition and Phenomenology of Space} edited by Luciano Boi and Carlos Lobo (Birkh\"auser, 2022).}
%\bigskip

\vspace{1cm}

\maketitle

\tableofcontents

\section{\bf  Introduction}
\label{sec:intro}

\begin{figure}%[h!] 
 %\vspace*{6mm}
 \centering 
 \includegraphics[scale=0.50]{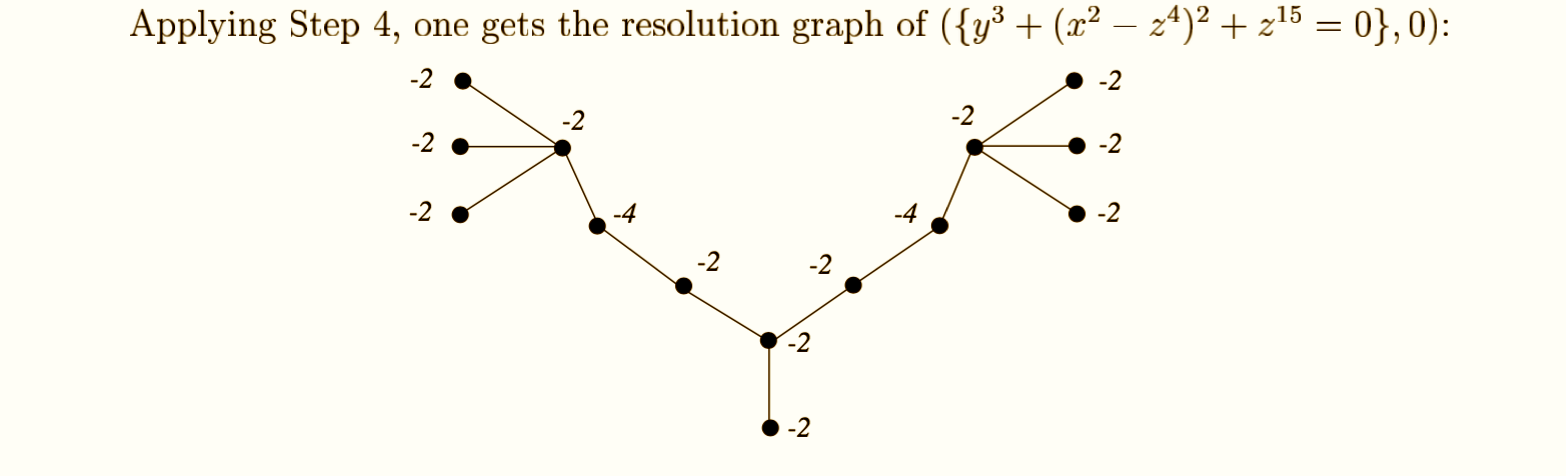} 
 \vspace*{1mm} 
 \caption{An example from a 2000 paper of N\'emethi and Szil\'ard} 
 \label{fig:2000-Nemethi}
 \end{figure}

\begin{figure}%[h!] 
% \vspace*{6mm}
 \centering 
 \includegraphics[scale=0.60]{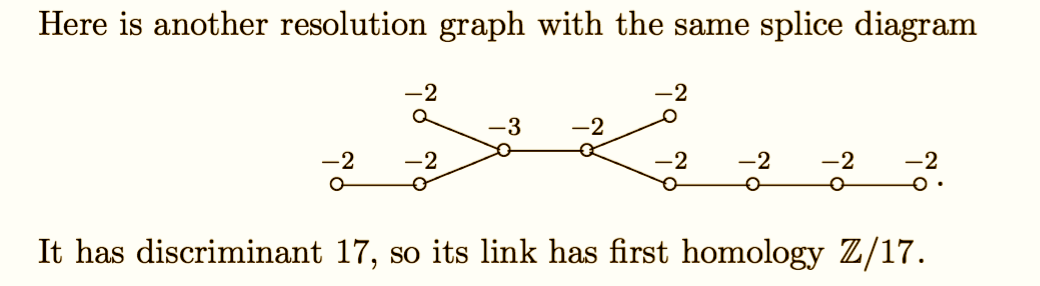} 
 \vspace*{1mm} 
 \caption{An example from a 2005 paper of Neumann and Wahl} 
 \label{fig:2005-Neumann-Wahl}
 \end{figure}

\begin{figure}%[h!] 
 %\vspace*{6mm}
 \centering 
 \includegraphics[scale=0.60]{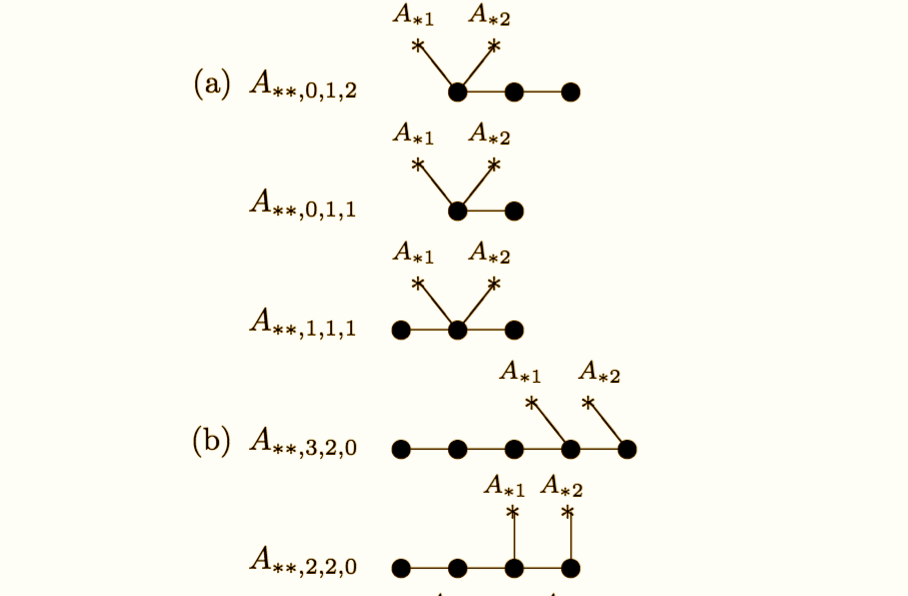} 
 \vspace*{1mm} 
 \caption{An example from a 2009 paper of Chung, Xu and Yau} 
 \label{fig:2009-Chung-Xu-Yau}
 \end{figure}

Nowadays, graphs are common tools in singularity theory. 
They mainly serve to represent morphological aspects of algebraic surfaces in the 
neighborhoods of their singular points.  
Three examples of such graphs may be seen in Figures 
\ref{fig:2000-Nemethi}, \ref{fig:2005-Neumann-Wahl}, \ref{fig:2009-Chung-Xu-Yau}.  
They are extracted from the papers \cite{NS 00}, \cite{NW 05} 
and \cite{CXY 09}, respectively. 

Comparing those figures, one may see that the vertices are diversely depicted by small stars or 
by little circles, which are either full or empty. These drawing conventions are not 
 important. What matters is that all vertices are decorated with numbers. I will explain 
 their meaning later. 

My aim in this paper is to understand which kinds of problems forced 
mathematicians to associate graphs to surface singularities. I will show how   
the initial idea, appeared in the 1930s, was only described in words, 
without any visual representation. 
Then I will suggest two causes that made the drawing of such  
graphs unavoidable, starting from the beginning of the 1960s. One 
of them was a topological reinterpretation of those graphs.  The other one was the 
growing interest in problems of classification of special types of surface singularities. 

Let me describe briefly the structure of the paper. In Section \ref{sec:meangraph}, I  
explain the meaning of such graphs: they represent configurations 
of curves which appear by resolving surface singularities. I continue in 
Section \ref{sec:wmres}  by describing what it means to ``resolve''  a singularity. 
In Section \ref {sec:resing} I present several models of surface singularities 
made around 1900 and I discuss one of the 
oldest configurations of curves, perhaps the most famous of them all: 
the 27 lines lying on a smooth cubic surface. In Section \ref{sec:DVCox}, I present 
excerpts of Du Val's 1934 paper in which he described a way of thinking  
about a special class of surface singularities in terms of  graphs. In the same paper,  
he made an analogy between his configurations of curves 
and the facets of special spherical simplices analyzed by Coxeter in his 1931 study 
of finite groups generated by reflections. It is perhaps the fact that Coxeter had described 
a way to associate a graph to such a simplex which, through this analogy, gave birth 
to Du Val's idea of speaking about graphs of curves. In Section \ref{sec:Mumpap}, 
I jump to the years 1960s, because until then Du Val's idea of associating 
graphs to singularities had almost never been used. I show how things changed with a
1961 paper of Mumford, in which he reinterpreted 
those graphs in the realm of $3$-dimensional topology.  
Hirzebruch's 1963 Bourbaki Seminar talk about this work of Mumford seems 
 to be the first place in which graphs representing arbitrary configurations of curves 
 were explicitly defined. I begin Section \ref{sec:WN} 
 by discussing  a 1967 paper of Waldhausen, 
 in which he built a subtle theory of the $3$-dimensional manifolds associated 
 to graphs as explained in Mumford's paper. I finish it with a discussion of a 1981 paper 
 of Neumann, which turned Waldhausen's work into a concrete ``calculus'' for deciding 
 whether two graphs represent the same $3$-dimensional manifold. 
 In Section \ref{sec:concl}, I conclude  
 by mentioning several recent directions of research concerning graphs associated to  
 singularities of algebraic varieties, and by summarizing this paper.

\medskip
{\bf Acknowledgements.} I am grateful to Luciano Boi, Franck Jedrzejewski 
and Carlos Lobos for the 
invitation to give a talk at the international conference ``\emph{Quand la forme devient 
substance : puissance des gestes, intuition diagrammatique et ph\'enom\'enologie de l'espace}'', 
which took place at Lyc\'ee Henri IV in Paris from 25 to 27 January 2018. This paper is an  
expanded version of my talk. I am also grateful to David Mumford for answering  
my questions about the evolution of the notion of dual graph and to 
Octave Curmi, Michael L\"onne and Bernard Teissier for their remarks.  
Special thanks are due to Mar\'{\i}a Angelica Cueto, Silvia De Toffoli and Fran\c{c}ois L\^e 
for their careful readings 
of a previous version of this paper and for their suggestions.

\section{\bf What is the meaning of such graphs?}
\label{sec:meangraph}

Let me begin by explaining the meaning of the graphs associated to singularities of surfaces. 
In fact, the construction is not specific to singularities, one may perform it 
whenever is given a configuration of curves on a surface. The rule is very simple:
  \begin{itemize}
       \item each curve of the configuration is represented by a vertex; 
       \item two vertices are joined by an edge whenever the corresponding 
          curves intersect.
  \end{itemize}
  
  A variant of the construction introduces as many edges between two vertices  
  as there are points in common between the corresponding curves. 
  
  Note that this construction reverses the dimensions of the input objects. Indeed, 
  the curves, which have dimension one, are represented by vertices of the graph, which 
  have dimension zero. Conversely, the intersection points of two curves, which 
  have dimension zero, are represented by edges of the graph, which have dimension one. 
  Remark also that an intersection point lies on a curve of the configuration if 
  and only if its associated edge of the graph contains the vertex representing the curve. 
  It is customary nowadays in mathematics to speak about ``duality'' whenever one has such 
  a dimension-reversing and inclusion-reversing correspondence 
  between parts of two geometric configurations. For this reason, one speaks here 
  about the ``\emph{dual graph}'' of the curve configuration, a habit which became common  
  at the end of the 1960s. 
  
  \begin{figure}%[h!] 
 %\vspace*{6mm}
 \centering 
 \includegraphics[scale=0.60]{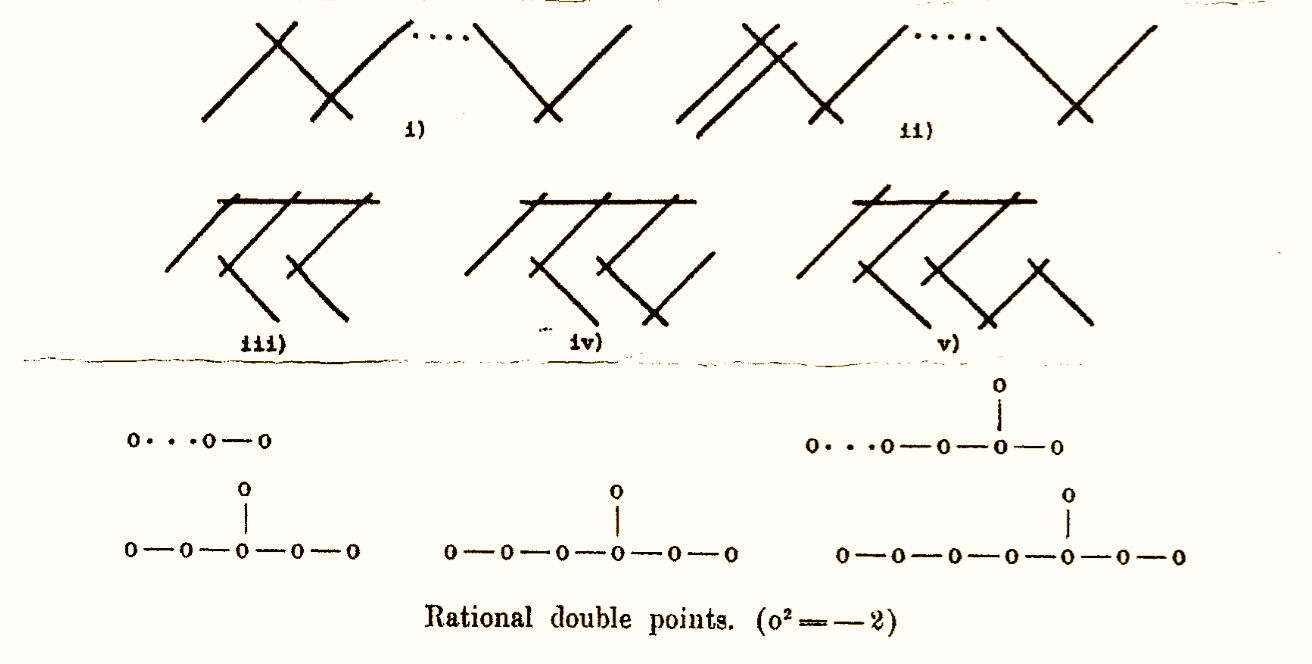} 
 \vspace*{1mm} 
 \caption{Artin's depictions of curve configurations and associated dual graphs} 
 \label{fig:1962-6-Artin}
 \end{figure}

  An example of this construction is represented in Figure \ref{fig:1962-6-Artin}, which 
  combines drawings from Michael Artin's 1962 and 1966 papers \cite{A 62, A 66}. 
  In the upper half, one sees sketches of curve configurations, each 
  curve being depicted as a segment. 
  This representation is schematic, as it does not respect completely the topology of 
  the initial curve configuration, which consists of curves without boundary points. 
  But it represents faithfully the intersections between the curves of the configuration: 
  two of its curves intersect if and only if the associated segments do. 
  The corresponding ``dual graphs'' are depicted in the lower part of the figure. For instance, 
  the vertex which is joined to three other vertices in the graph of the lower right corner represents 
  the horizontal segment of the curve configuration labeled v), on the right of the second row. 
  
  Note that the representation of the curve configurations as dual graphs emphasizes better 
  visually its overall connectivity pattern than the representation as a configuration of segments. 
  This is probably one of the reasons which led Artin to pass from drawings of 
  configurations of segments  in his 1962 paper to drawings of dual graphs in his 1966 paper. 

 More generally, dual graphs may be introduced whenever one is interested 
 in the mutual intersections of several subsets of a given set. It is not important that 
 the given sets consist of the points of several curves lying on surfaces, they may for instance be 
arbitrary subsets of manifolds of any dimension or, less geometrically, the sets of members 
 of various associations of persons. Then, one represents each set by a vertex  
 and one joins two such vertices by an edge if and only if the corresponding sets 
 intersect. 
 
 As a general rule, one represents  any object of study by a vertex,  
 whenever one is not interested in its internal structure, but in its 
 ``sociology'', that is, in its relations or interactions with other objects. A basic way to depict  
 these relations is to join two such vertices whenever the objects represented by them 
 interact in the way under scrutiny. In our context, one considers that two curves interact 
 if and only if they have common points. 
 
 Let us come back to the configurations of curves and their associated dual graphs from  
 Figure \ref{fig:1962-6-Artin}. Those drawings represent the classification of a special type  
 of surface singularities, namely the ``\emph{rational double points}'', in the terminology of 
 Artin's 1966 paper \cite{A 66}. That list was not new, it had already appeared 
 in the solution of a different classification problem -- leading nevertheless to the same objects -- 
 in Du Val's 1934 paper in which he had informally introduced  the idea 
 of dual graph. We will further discuss that paper in Section \ref{sec:DVCox}. 
 
 \begin{figure}%[h!] 
 %\vspace*{6mm}
 \centering 
 \includegraphics[scale=0.60]{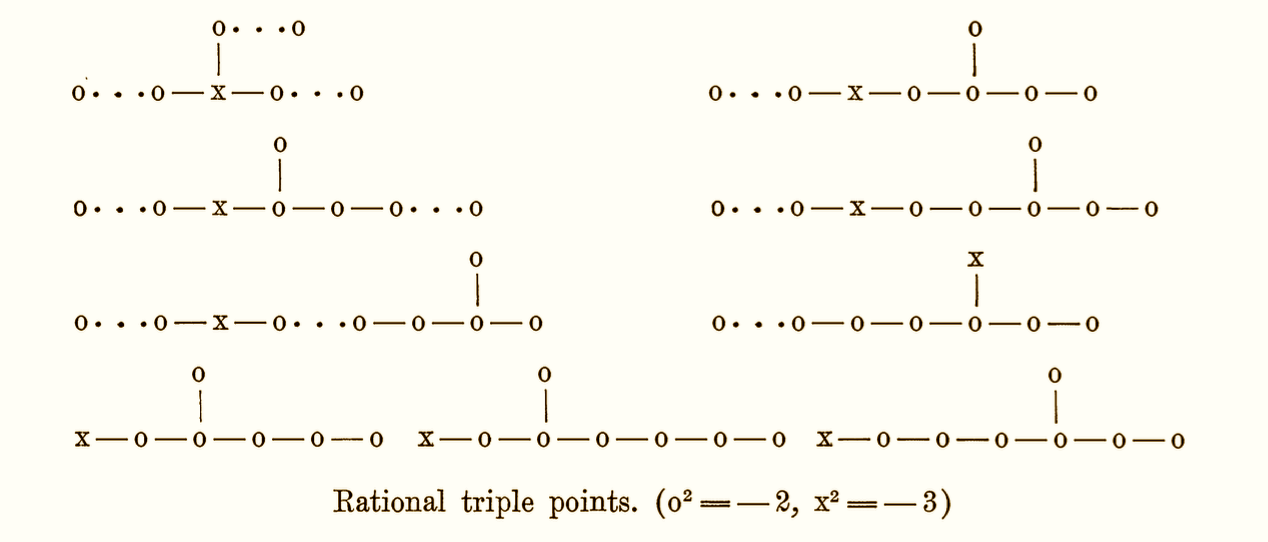} 
 \vspace*{1mm} 
 \caption{Artin's classification of dual graphs of rational triple points} 
 \label{fig:1966-Artin-rtp}
 \end{figure}

 Before passing to the next section, 
 let me mention that Artin's paper \cite{A 66} contained also a new classification, that of the 
 dual graphs associated to the ``\emph{rational triple points}'' (see Figure \ref{fig:1966-Artin-rtp}). 
 Being more abundant 
 than those of the lower part of Figure \ref{fig:1962-6-Artin}, it becomes apparent 
 that it is also more economical for printing to draw such graphs rather 
 than configurations of segments. 
 
 Now that we understood how configurations of curves lead to dual graphs, let us see in which way 
 singularities of surfaces may lead to configurations of curves. This is the object of the next section.

\section{\bf What does it mean to resolve the singularities of an algebraic surface?}
\label{sec:wmres}

\begin{figure}%[h!] 
 %\vspace*{6mm}
 \centering 
 \includegraphics[scale=0.90]{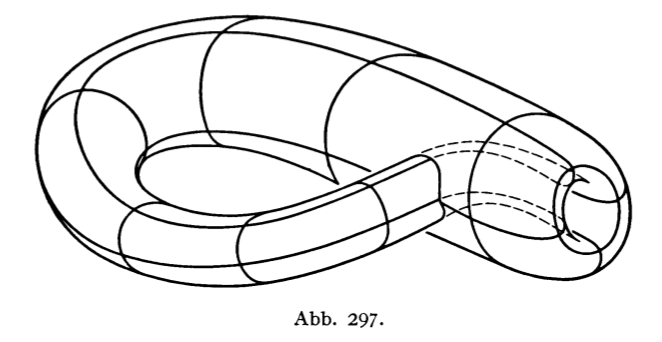} 
 \vspace*{1mm} 
 \caption{A Klein bottle} 
 \label{fig:1932-Klein-bottle}
 \end{figure}

\begin{figure}%[h!] 
 \vspace*{6mm}
 \centering 
 \includegraphics[scale=0.90]{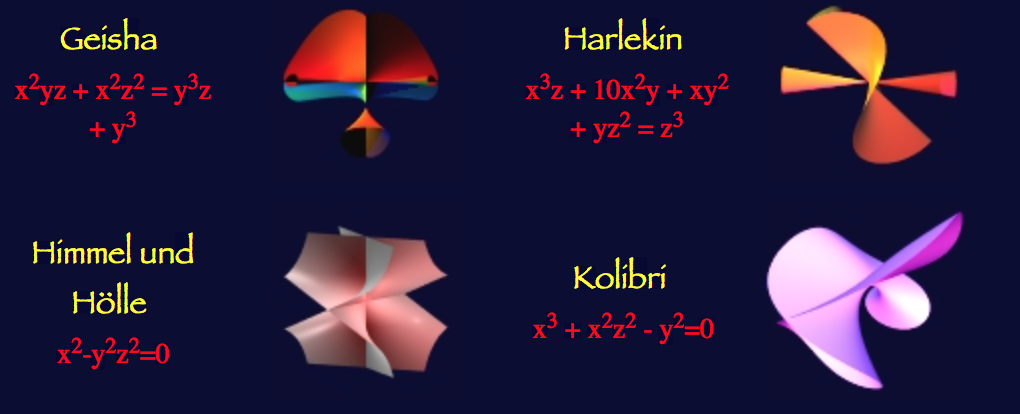} 
 \vspace*{1mm} 
 \caption{Several real surface singularities from Hauser's gallery} 
 \label{fig:2018-Hauser}
 \end{figure}

What is a \emph{singularity of an algebraic surface}? It is a special point, 
at which the surface is not smooth. 
For instance, a sphere does not have singular points, but a double cone, idealization of the boundary of 
a nighty region illuminated by a lighthouse, has a singular point at its vertex. In this case, 
the singular point is isolated, but other surfaces may have whole curves of singularities.  
Such curves may be either self-intersections of the surface, 
as shown in Figure \ref{fig:1932-Klein-bottle}\footnote{This illustration of an immersion of a Klein bottle 
 in three-dimensional cartesian space comes from the wonderful  
  1932 book \cite{HC 32} of David Hilbert and Stefan Cohn-Vossen. Note that in 
  addition to a circle of self-intersection, this illustration represents also curves of 
  apparent contours and several closed curves drawn on the surface.}, or 
they may exhibit more complicated behaviour, 
as shown in Figure \ref{fig:2018-Hauser}\footnote{This figure comes from \href{https://homepage.univie.ac.at/herwig.hauser/gallery.html}{https://homepage.univie.ac.at/herwig.hauser/gallery.html}, January 2018. Copyright: Herwig Hauser, University of Vienna, www.hh.hauser.cc.}. 

In this last figure,  a polynomial equation in three variables is written next to each surface.
The reason is that each of those surfaces is a portion of the locus of points which satisfy 
the associated equation 
in the $3$-dimensional cartesian space of coordinates $(x,y,z)$. As polynomials 
are algebraic objects, such a locus is called ``\emph{algebraic}''. 
There are also algebraic surfaces in cartesian spaces of higher 
dimensions, defined by systems of polynomial equations in more than three variables. 
In fact, all surfaces considered in the papers discussed here are algebraic. One advantage of 
working with such surfaces is that one may consider not only 
the real solutions of those equations, but also the complex ones. 
In this way, one expects in general to make the correspondence between 
the algebraic properties of the defining equations and the morphological properties of the 
associated surface easier to understand. 

A prototype of this expectation is the fact that a polynomial equation 
in one variable has as many complex roots as its degree, provided that the roots are counted with 
suitable ``multiplicities'' (this is the so-called ``fundamental 
theorem of algebra'', but it is rather a fundamental theorem of the \emph{correspondence} between 
algebra and topology). If one considers instead only its real roots, then their number is 
not determined by the degree, but there are several possibilities. In fact, as I will briefly 
explain at the beginning of Section \ref{sec:resing}, whenever one considers 
families with three parameters of polynomials, these possibilities may be distinguished using 
``\emph{discriminant surfaces}'', which have in general non-empty singular loci. 

Because of this expected relative simplicity of complex algebraic geometry versus real algebraic 
geometry,  it became customary in the XIXth century to study the sets of \emph{complex} solutions 
of polynomial equations in three variables. One gets in this way ``\emph{complex algebraic surfaces}''. 
Nevertheless, in order to build an intuition of their properties, it may be useful to 
practice with concrete \emph{models} of the associated \emph{real} surfaces. 
Around the end of the XIXth century, such models were either drawn or 
manufactured using for instance wood, plaster, cardboard, wires and string. Nowadays they are 
also built using 3D-printers or, more commonly, simulated using techniques of computer visualization. 
This is for instance the case of Hauser's images of Figure \ref{fig:2018-Hauser}.

Why is it important to study singular surfaces? Because, in general, surfaces do not appear  
alone, but rather in families depending on parameters (which, in physical contexts, may be 
for instance temperatures 
or intensities of external fields), and that for some special values of these parameters one gets 
surfaces with singularities. Understanding the singular members of a family is 
many times essential for understanding also subtle aspects of its non-singular members. 
For instance, one may understand part of the structure of a non-singular member by looking 
at its portions which ``vanish'' when one converges to a singular member.

The techniques of differential or algebraic geometry used in the study of 
smooth algebraic surfaces may be extended to singular surfaces using three basic procedures:
   \begin{itemize} 
         \item by decomposing a singular surface into smooth ``\emph{strata}'', which are either 
           isolated points, smooth portions of curves or smooth pieces of surfaces; this is 
           similar to the decomposition of the surface of a convex polyhedron into vertices, 
           edges and faces; 
         \item by seeing a singular surface as a limit of smooth ones; when this is possible, 
             one says that the surface was ``\emph{smoothed}''; such a process is not always 
             possible, and even if it is possible, it can be usually done in various ways; 
         \item by seeing a singular surface as a projection of a smooth one, living 
              in a higher dimensional ambient space; if such a projection 
              leaves the smooth part of the initial singular surface unchanged, then 
              it is called a ``\emph{resolution of singularities}''; resolutions of singularities always 
              exist, but are not unique.   
   \end{itemize} 
   
   Intuitively speaking, resolving the singularities of a surface means to remove its singular locus and 
   to replace it algebraically by another configuration of points and curves, so that the 
   resulting surface is smooth. In the special case of 
   an \emph{isolated} singular point of algebraic surface,  
   one replaces that singular point by a configuration of curves, called the 
   ``\emph{exceptional divisor}'' of the resolution. For instance, all the graphs 
   appearing in Figures \ref{fig:2000-Nemethi}--\ref{fig:1966-Artin-rtp}  are dual graphs of 
   exceptional divisors of resolutions of isolated singularities of complex surfaces. 
   
   In simple examples, one may resolve the singularities of an algebraic surface 
   by performing finitely many times 
   the elementary operation called ``\emph{blowing up a point}'', which is a mathematical way to 
   look through a microscope at the neighborhood of the chosen point. 
   This operation builds new cartesian 
   spaces starting from the space which contains the initial singular surface. 
   Each of those spaces contains a new surface, which projects onto a part of the initial one. If one 
   of those surfaces is smooth, then one keeps it untouched. Otherwise, one blows up again its 
   isolated singular points. It may happen that finitely many such operations lead to a  
   family of 
   {\em smooth} surfaces, each one of them projecting onto a portion of the initial singular surface. 
   Those surfaces may 
   be glued, together with their projections, into a global smooth surface which ``\emph{resolves 
   the singularities}'' of the initial one. 
   
 \begin{figure}%[h!] 
 %\vspace*{6mm}
 \centering 
 \includegraphics[scale=0.80]{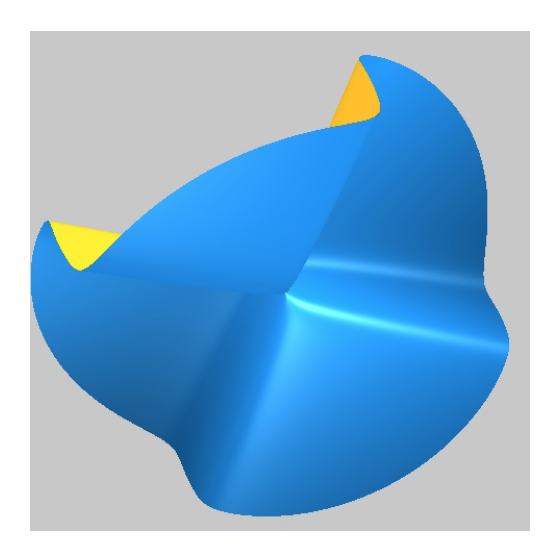} 
 \vspace*{1mm} 
 \caption{An isolated $A_6$ surface singularity} 
 \label{fig:2018-Fruhbis-Kruger-1}
 \end{figure}

\begin{figure}%[h!] 
 %\vspace*{6mm}
 \centering 
 \includegraphics[scale=0.40]{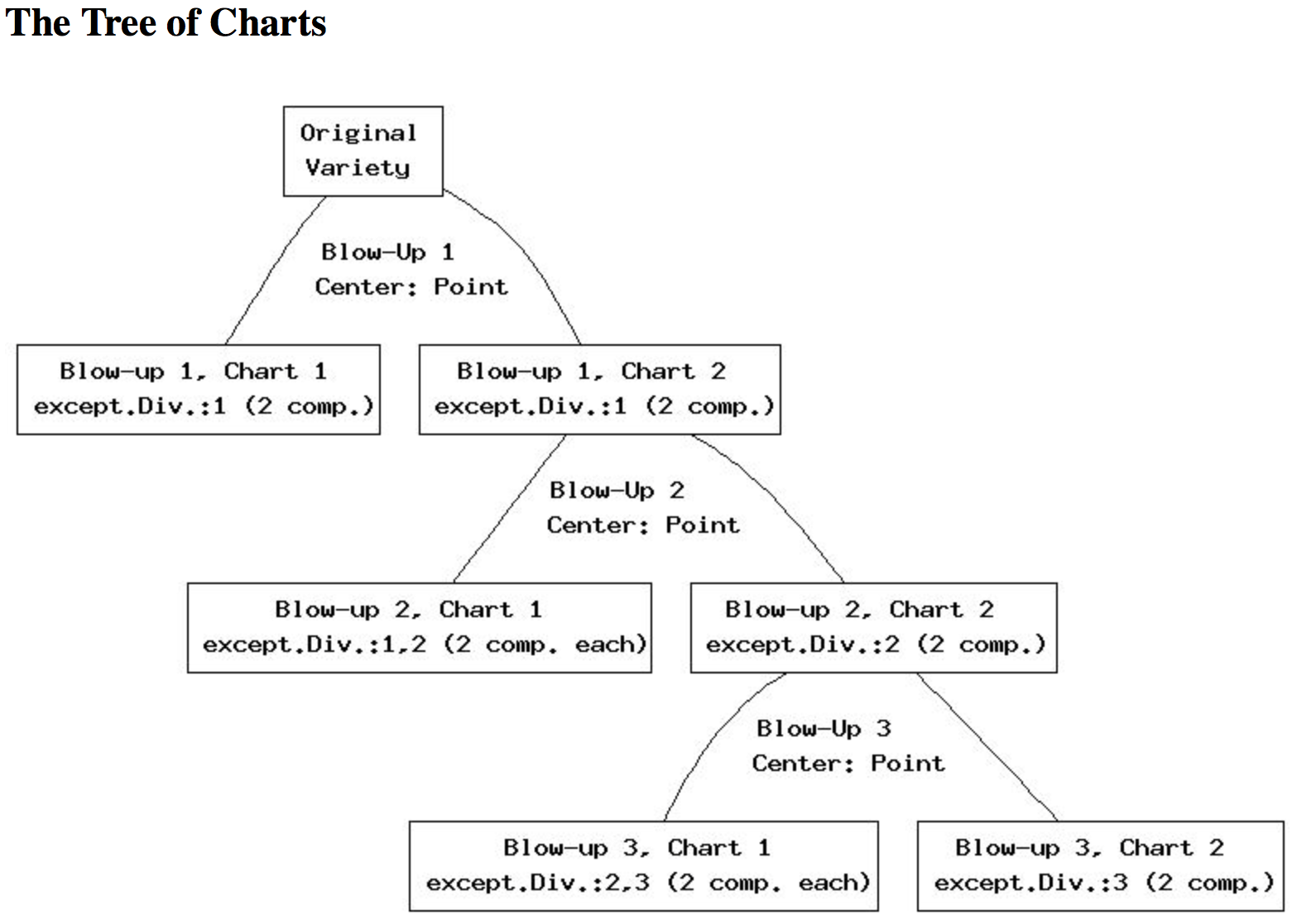} 
 \vspace*{1mm} 
 \caption{Fr\"uhbis-Kr\"uger's representation of a resolution process} 
 \label{fig:2018-Fruhbis-Kruger-2}
 \end{figure}

\begin{figure}%[h!] 
 %\vspace*{6mm}
 \centering 
 \includegraphics[scale=0.50]{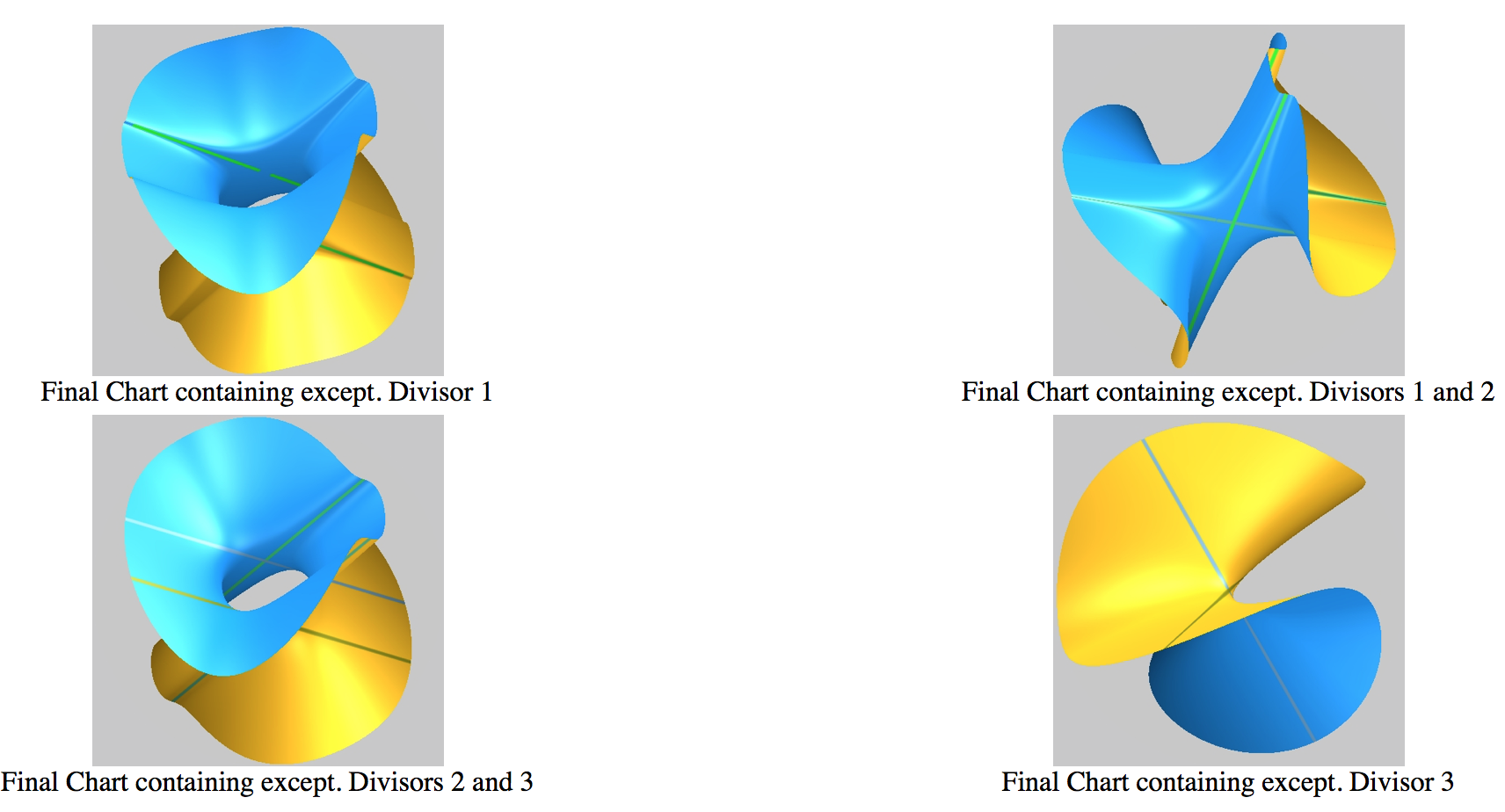} 
 \vspace*{1mm} 
 \caption{The final stages of the previous resolution process} 
 \label{fig:2018-Fruhbis-Kruger-3}
 \end{figure}
   
   Let us consider for instance the surface with equation $x^2 - y^2 + z^7 =0$, illustrated 
   in Figure \ref{fig:2018-Fruhbis-Kruger-1}. 
   It has an isolated singularity (called ``of type $A_6$'') at the 
   origin\footnote{This figure was taken  
   from the web-page https://www.krueger-berg.de/anne/aufl-bilder/A6.html   of 
   Anne Fr\"uhbis-Kr\"uger in September 2020. This is also the case of Figures 
   \ref{fig:2018-Fruhbis-Kruger-2} and \ref{fig:2018-Fruhbis-Kruger-3}. }. 
    If one performs the previous iterative process of blowing up the 
    singular points of the intermediate surfaces, one gets a ``tree'' of surfaces, 
    represented diagrammatically on Figure \ref{fig:2018-Fruhbis-Kruger-2}. The initial surface is 
    indicated in the top-most rectangle, and each edge of the diagram represents a 
    blow-up operation.

The final smooth surfaces produced by the process are represented 
in Figure \ref{fig:2018-Fruhbis-Kruger-3}.
Each of them contains one or more highlighted lines. Those lines glue into a configuration of 
curves on the total smooth surface which resolves the initial singular one. 
This configuration is the exceptional divisor of this resolution of the 
starting isolated singularity. By looking carefully at the way the gluing is performed, one 
may show that its associated dual graph is a chain of five segments. This 
means that it is of the type shown on the left of the third row in Figure \ref{fig:1962-6-Artin}.

 For more complicated singularities, it may not be enough to blow up points, as previous 
 blow-ups may create whole curves of singularities. Other operations which allow to 
 modify the singular locus were introduced in order to deal with this problem. One 
 may learn about them in Koll\'ar's book \cite{K 07}, which 
 explains various techniques of resolution of singularities in any dimension. The reader 
 more interested in gaining intuition about resolutions of surfaces may consult 
 Faber and Hauser's promenade \cite{FH 10} through a garden of examples of resolutions 
 or my introduction \cite{PP 11bis} to one of the oldest methods of resolution, 
 originating in Jung's method for parametrizing algebraic surfaces locally.

\section{\bf Representations of surface singularities around 1900}
\label{sec:resing}

In the previous sections we saw contemporary representations of surface singularities, obtained 
using computer visualization techniques. Let us turn now to older representations, dating 
back to the beginning of the XXth century. This will allow us to present one of the oldest sources 
of surfaces with singularities -- the study of discriminant surfaces -- and one of the most 
famous configurations of curves -- consisting of the lines contained in a smooth cubic surface. 

\begin{figure}%[h!] 
 %\vspace*{6mm}
 \centering 
 \includegraphics[scale=0.6]{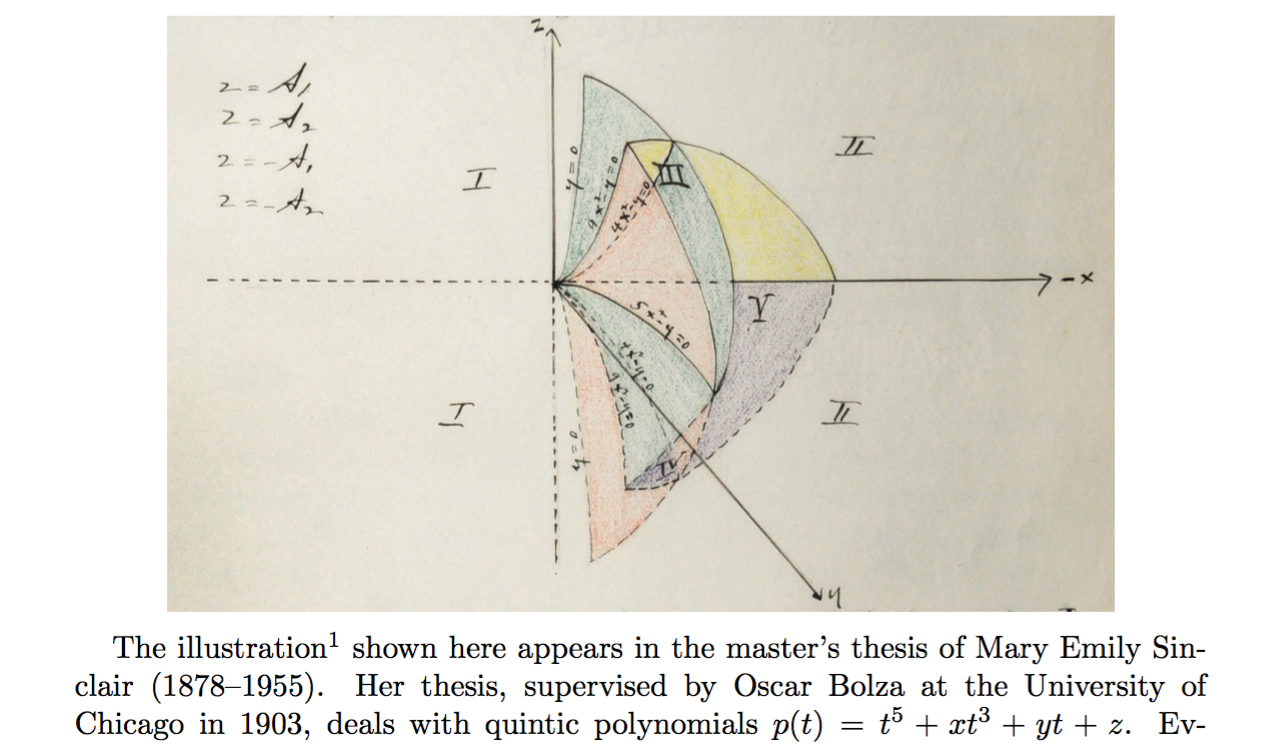} 
 \vspace*{1mm} 
 \caption{Sinclair's representation of a discriminant surface} 
 \label{fig:1903-Sinclair}
 \end{figure}

Figure \ref{fig:1903-Sinclair} shows\footnote{This drawing and the text 
immediately below it were extracted from the paper \cite{TW 11} of Jaap Top and Erik Weitenberg.} 
a hand-drawn ``\emph{discriminant surface}'', which 
has whole curves of singular points, as was the case in the examples of Figures 
\ref{fig:1932-Klein-bottle} and \ref{fig:2018-Hauser}. 
It reproduces a drawing done by Mary Emily Sinclair in her 1903 thesis. Let me discuss 
this surface a little bit, as it emphasizes another source of interest on the structure of singular 
surfaces. As explained in its caption from the paper \cite{TW 11}, Sinclair was studying the 
family with three parameters $(x,y,z)$ of polynomials of the form $t^5 + xt^3 + yt + z$. 
One may associate to it an 
algebraic family of sets of points, namely the sets of roots of the polynomial in the variable $t$ 
obtained for fixed  values of the parameters. 
The ``discriminant surface'' is the subset of the cartesian space of coordinates 
$(x,y,z)$ for which the associated polynomial has at least one multiple complex root. 

More generally, consider any 
family of points, curves, surfaces or higher-dimensional algebraic objects, depending 
algebraically on some parameters. If there are exactly three parameters, then the set of 
singular objects of the family is usually a surface in the space of 
parameters. All the surfaces obtained in this way are called ``{\em discriminant surfaces}'', because 
they allow to discriminate the possible aspects of the objects in the family, according 
to the position of the corresponding point in the space of parameters, relatively to the surface. 
For instance, by determining in which region of the complement of the surface of 
Figure \ref{fig:1903-Sinclair} lies the point with coordinates $(x,y,z)$, one may see if the set of real roots 
of the polynomial has $1, 3$ or $5$ elements -- those 
being the only possibilities for a quintic polynomial equation, because the non-real 
roots come in pairs of complex conjugate numbers. The reader interested in the analogous 
study of quartic polynomial equations may read Michel Coste's paper \cite{C 10}.

  \begin{figure}%[h!] 
 %\vspace*{6mm}
 \centering 
 \includegraphics[scale=0.6]{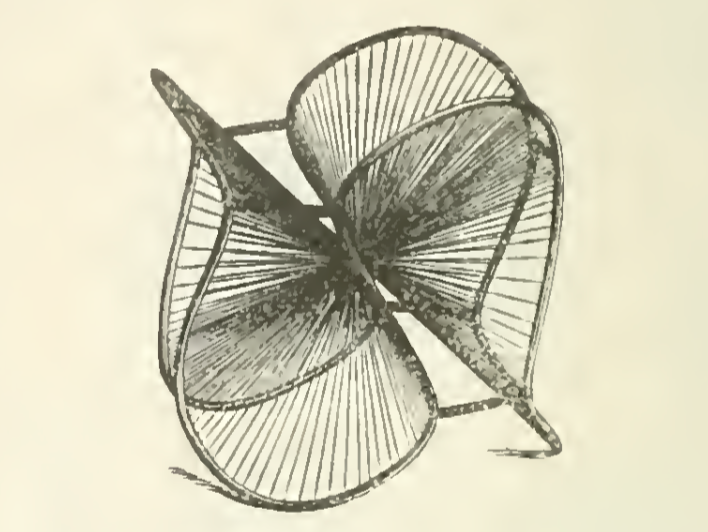} 
 \vspace*{1mm} 
 \caption{An image from Schilling's catalog of mathematical models} 
 \label{fig:1903-Schilling}
 \end{figure}
 
\begin{figure}%[h!] 
 %\vspace*{6mm}
 \centering 
 \includegraphics[scale=0.55]{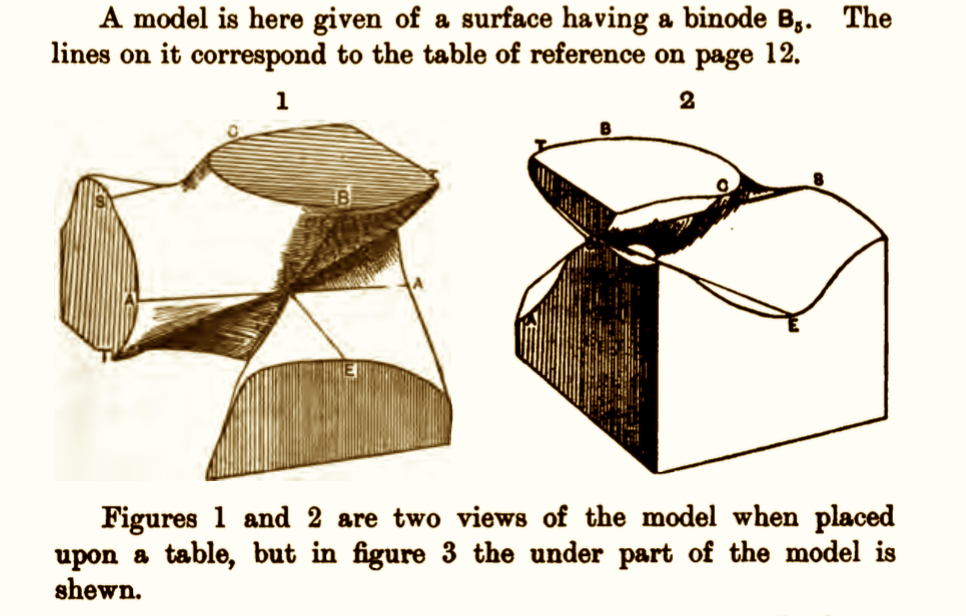} 
 \vspace*{1mm} 
 \caption{An image from Blythe's book on models of cubic surfaces} 
 \label{fig:1905-Blythe}
 \end{figure}

 Let us pass now to material models of surfaces with singularities. Figure 
\ref{fig:1903-Schilling} reproduces an engraving\footnote{It may be found on page 123, 
part II.3.b of the catalog \cite{* 11}.} from the 1911 catalog of mathematical models of Martin Schilling's 
enterprise. It depicts a cone over a smooth cubic curve, that is, a smooth curve contained 
in the projective plane and defined by the vanishing of a homogeneous polynomial in $3$ variables. 
This cone has 
therefore only one singular point, its vertex. Figure \ref{fig:1905-Blythe} 
shows a reproduction from the 1905 book \cite{B 05} of William Henry Blythe. 
It depicts two plaster models of cubic surfaces with singularities.

 \begin{figure}%[h!] 
 %\vspace*{6mm}
 \centering 
 \includegraphics[scale=0.5]{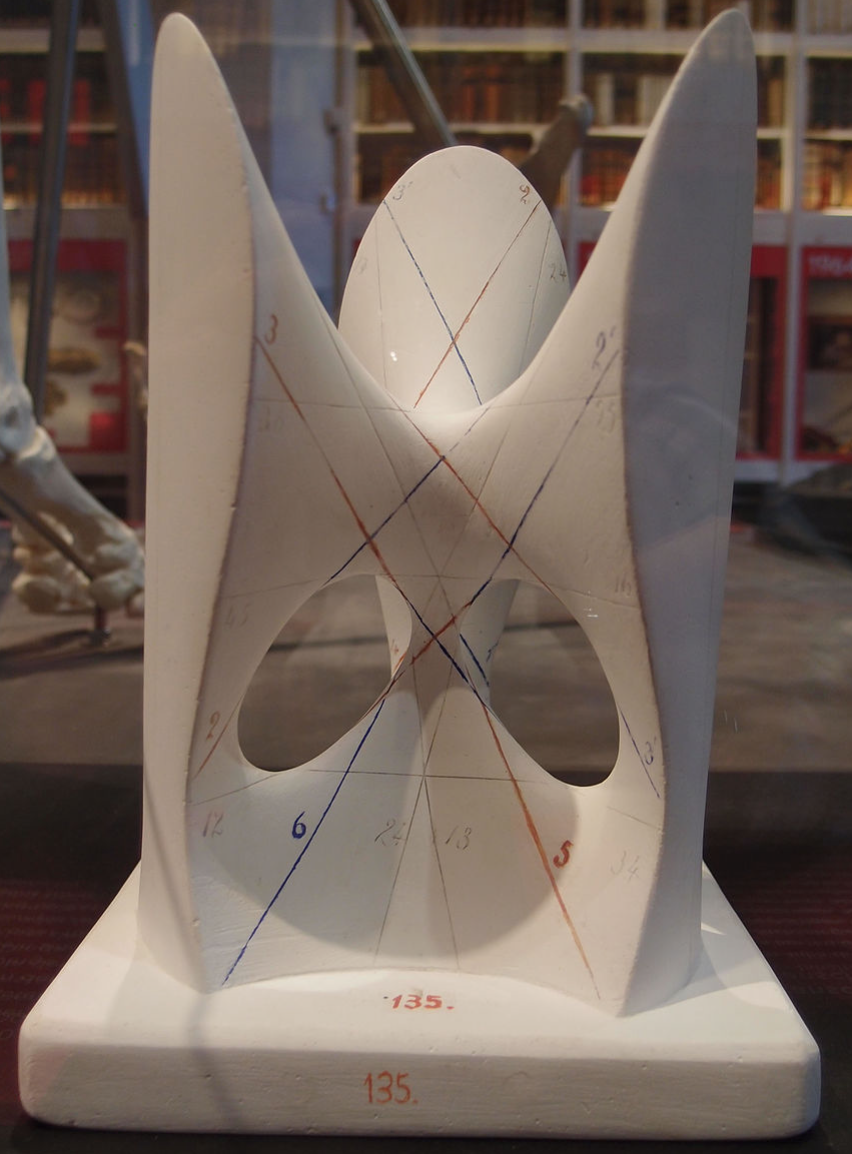} 
 \vspace*{1mm} 
 \caption{A model of Clebsch's diagonal surface and of its $27$ lines} 
 \label{fig:2012-Zauzig-Clebsch-diag}
 \end{figure}
 
One of the most famous discoveries of the XIXth century regarding the properties 
of algebraic surfaces is that all smooth complex algebraic cubic surfaces 
situated in the projective space of dimension three contain exactly 27 lines. This discovery 
was done in 1849, during a correspondence between Arthur Cayley and George Salmon, and 
it triggered a lot of research\footnote{For instance, in the historical 
summary of his 1915 thesis \cite{H 15}, Henderson mentions that ``in a bibliography 
on curves and surfaces compiled by J. E. Hill [in 1897] [...] the section on cubic surfaces contained 
two hundred and five titles. The Royal Society of London Catalogue of Scientific Papers, 
1800--1900, volume for \emph{Pure Mathematics} (1908), contains very many more.''}. 
Starting from around 1870, material models of parts of real cubic surfaces with all 
$27$ lines visible on them started to be built. One may see such a model 
in Figure \ref{fig:2012-Zauzig-Clebsch-diag}. It represents a portion of 
``\emph{Clebsch's diagonal surface}''\footnote{Clebsch's diagonal surface is usually defined by the 
pair of homogeneous equations $x_0 + \cdots + x_4 = 0$ and $x_0^3 + \cdots + x_4^3 = 0$ inside the 
projective space of dimension $4$ whose homogeneous coordinates are denoted 
$[x_0 : \cdots : x_4]$. This photograph of a model belonging to the 
University of G\"ottingen was taken by Zausig in 2012. It comes from Wikimedia Commons: 
https://commons.wikimedia.org/wiki/File:Modell\_der\_Diagonalfläche\_von\_Clebsch\_-Schilling\_VII,\_1\_-\_44-.jpg}.  This surface does not contain singular points, but it is interesting 
in our context because it 
 exhibits a highly sophisticated configuration of curves, composed of its $27$ lines. 

Much more details about the building of models of algebraic surfaces around 1900 may 
be found in the books \cite{* 17} and \cite{* 17bis}. As illustrated by  Figure \ref{fig:1905-Blythe}, 
plaster models 
were built not only of smooth cubic surfaces, but also of singular ones. The 
manufacturing process was based on Rodenberg's 1878 work \cite{R 78}. The complete 
classification of the topological types of real cubic surfaces was achieved by 
Kn\"orrer and Miller in their 1987 paper \cite{KM 87}. Other historical details about the study 
of the configurations of 27 lines lying on smooth cubic surfaces may be found in 
Polo Blanco' and L\^e's theses \cite{PB 07} and \cite{L 15}, as well as in L\^e's paper \cite{L 13} 
and Labs' paper \cite{L 17}. 

Note that at the beginning of the XXth century, 
some artists from the Constructivist and Surrealist movements were inspired 
by material models of possibly singular surfaces, as explained in Vierling-Claassen's article \cite{VC 17}. 
It would be interesting to know in which measure computer models as those of 
Figure \ref{fig:2018-Hauser} inspire nowadays other artists.

\section{\bf Du Val's singularities, Coxeter's diagrams and the birth of dual graphs}
\label{sec:DVCox}

Let us discuss now the 1934 paper \cite{V 34} in which Patrick Du Val considered, seemingly  
for the first time, the idea of \emph{dual graph} of an exceptional divisor of resolution 
of surface singularity. 

As indicated by the title of his paper, 
Du Val's problem was to classify the ``\emph{isolated singularities of surfaces which do not affect 
the conditions of adjunction}''. Given a possibly singular algebraic surface 
contained in a complex projective space of dimension three, its ``\emph{adjoint surfaces}'' are 
other algebraic surfaces contained in the same projective space and defined in terms of double integrals. 
I will not give here their precise definition, 
which is rather technical\footnote{The interested reader may find it in Merle and Teissier's paper \cite{MT 77}. The whole volume containing that paper is dedicated 
to a modern study of the singularities analyzed by Du Val.}. Let me only mention 
that the adjoint surfaces must contain all curves consisting of singularities of the given surface. 
By contrast, it does not necessarily contain its \emph{isolated} singular points. 
Those through which the adjoint 
surfaces are not forced to pass are precisely the singularities 
``which do not affect the conditions of adjunction''.

\begin{figure}%[h!] 
 %\vspace*{6mm}
 \centering 
 \includegraphics[scale=0.6]{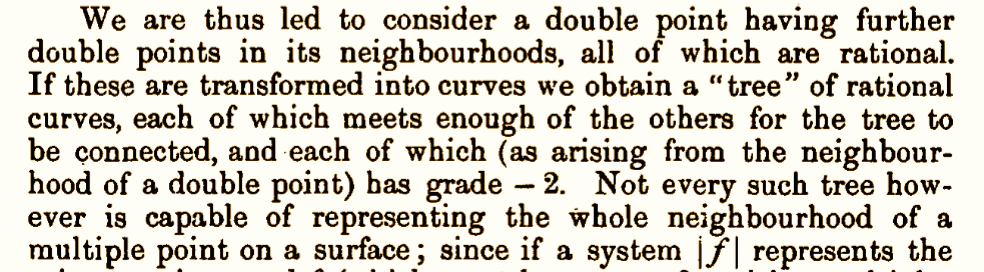} 
 \vspace*{1mm} 
 \caption{Du Val's introduction of dual graphs} 
 \label{fig:1934-Du-Val-graphes-duaux}
 \end{figure}

\begin{figure}%[h!] 
 %\vspace*{6mm}
 \centering 
 \includegraphics[scale=0.6]{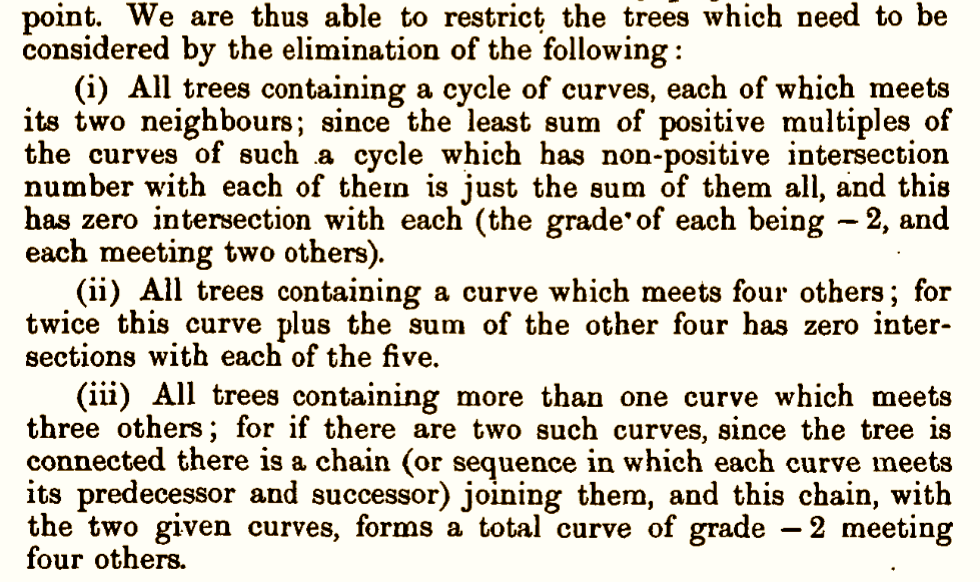} 
 \vspace*{1mm} 
 \caption{Du Val's restricted class of graphs} 
 \label{fig:1934-Du-Val-pas-encore}
 \end{figure}

 Du Val analyzed such singularities by looking at their resolutions. It is in this context that 
he wrote that for each one of those singularities, there is a resolution whose 
associated exceptional divisor is a ``\emph{``tree'' of rational curves}'' 
with supplementary properties (see Figure \ref{fig:1934-Du-Val-graphes-duaux}). For instance, 
each curve in this ``tree''  has necessarily self-intersection $-2$ in its ambient smooth surface 
(this is the meaning of the syntagm ``\emph{has grade $-2$}''). Du Val continued by giving a list 
of constraints verified by such ``trees'', if they were to correspond to singularities which 
do not affect the conditions of adjunction (see Figure \ref{fig:1934-Du-Val-pas-encore}). 
Using those constraints, he arrived exactly at the 
list of configurations of curves depicted in Figure \ref{fig:1962-6-Artin}. But, in contrast with 
Artin's papers \cite{A 62, A 66} from 1960's, his article does not contain any schematic drawing of 
a configuration of curves, or of an associated dual graph. 

 It is not even clear whether Du Val really thought about dual graphs. Perhaps he drew for 
 himself some diagrams resembling those of the upper part of Figure \ref{fig:1962-6-Artin}, 
 and he saw an analogy with some ``trees'' considered by other mathematicians. Note that 
 it is possible that for Du Val the term ``tree'' meant what we call ``graph''. Indeed, one sees him 
 stating in the excerpt of Figure \ref{fig:1934-Du-Val-pas-encore} that the ``trees'' under 
 scrutiny should not contain ``a cycle of curves'', a formulation which allows some ``trees of curves'' 
 to contain such cycles. 
 
 \begin{figure}%[h!] 
 %\vspace*{6mm}
 \centering 
 \includegraphics[scale=0.6]{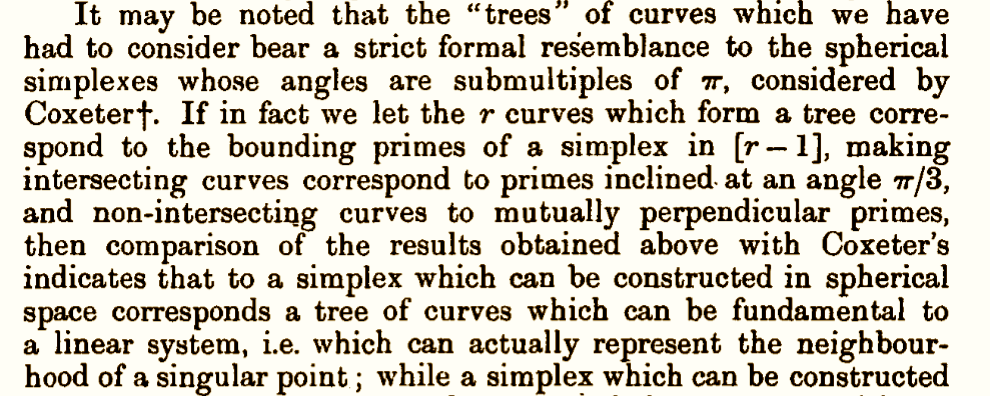} 
 \vspace*{1mm} 
 \caption{Du Val's analogy with Coxeter's spherical simplices} 
 \label{fig:1934-Du-Val-analogie}
 \end{figure}

At the end of his article, Du Val mentioned an analogy with results of  
Coxeter\footnote{Du Val cites the paper \cite{C 32}, but Coxeter already 
considered this problem one year earlier, in \cite{C 31}. Note that Du Val and 
Coxeter were friends and that they discussed regularly about their research. 
One may learn a few details 
about their friendship and discussions in Roberts' book 
\cite{R 07}, especially on pages 71-72.} regarding  
finite groups generated by reflections (see Figure \ref{fig:1934-Du-Val-analogie}). 
In order to understand this analogy, we have to know 
that Coxeter started from a finite set of hyperplanes passing through the 
origin in a real Euclidean vector space of arbitrary finite dimension. He assumed that they 
were spanned by the facets of a simplicial cone emanating from the origin, and he 
looked at the spherical simplex obtained by intersecting the cone with the unit sphere 
centered at the origin. Coxeter's problem was to classify those spherical simplices 
for which the group generated by the orthogonal reflections in the given hyperplanes 
is \emph{finite}. 

Du Val realized that his classification of isolated singularities which do not affect the conditions 
of adjunction corresponds to a part of Coxeter's classification of spherical simplices giving rise 
to finite groups of reflections.  In order to make this correspondence visible, 
he associated to each curve of a given exceptional 
divisor a facet of the simplex, two curves being disjoint if and only if the corresponding facets 
are orthogonal, and having one point of intersection if and only if the facets meet at an angle of 
$\pi/3$. 

\begin{figure}%[h!] 
 %\vspace*{6mm}
 \centering 
 \includegraphics[scale=0.55]{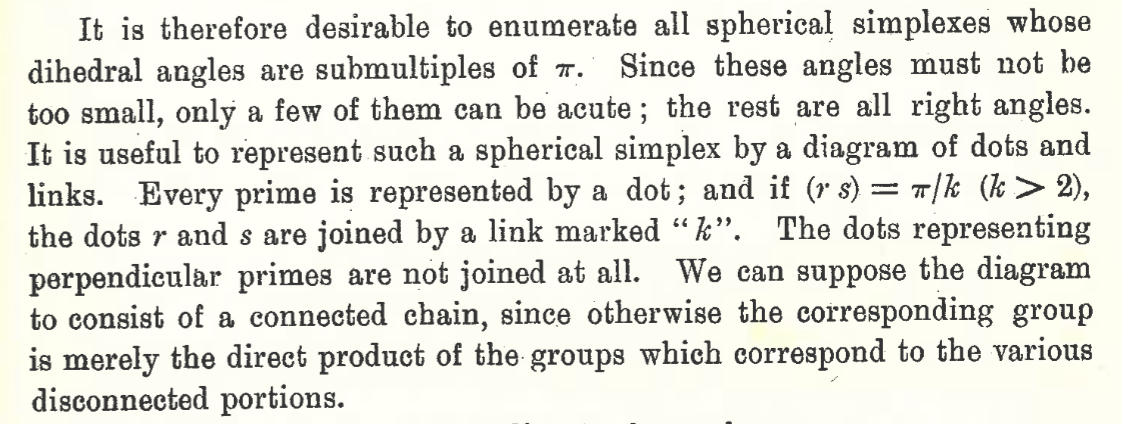} 
 \vspace*{1mm} 
 \caption{Coxeter's introduction of his diagrams} 
 \label{fig:1931-Coxeter}
 \end{figure}

\begin{figure}%[h!] 
 %\vspace*{6mm}
 \centering 
 \includegraphics[scale=0.5]{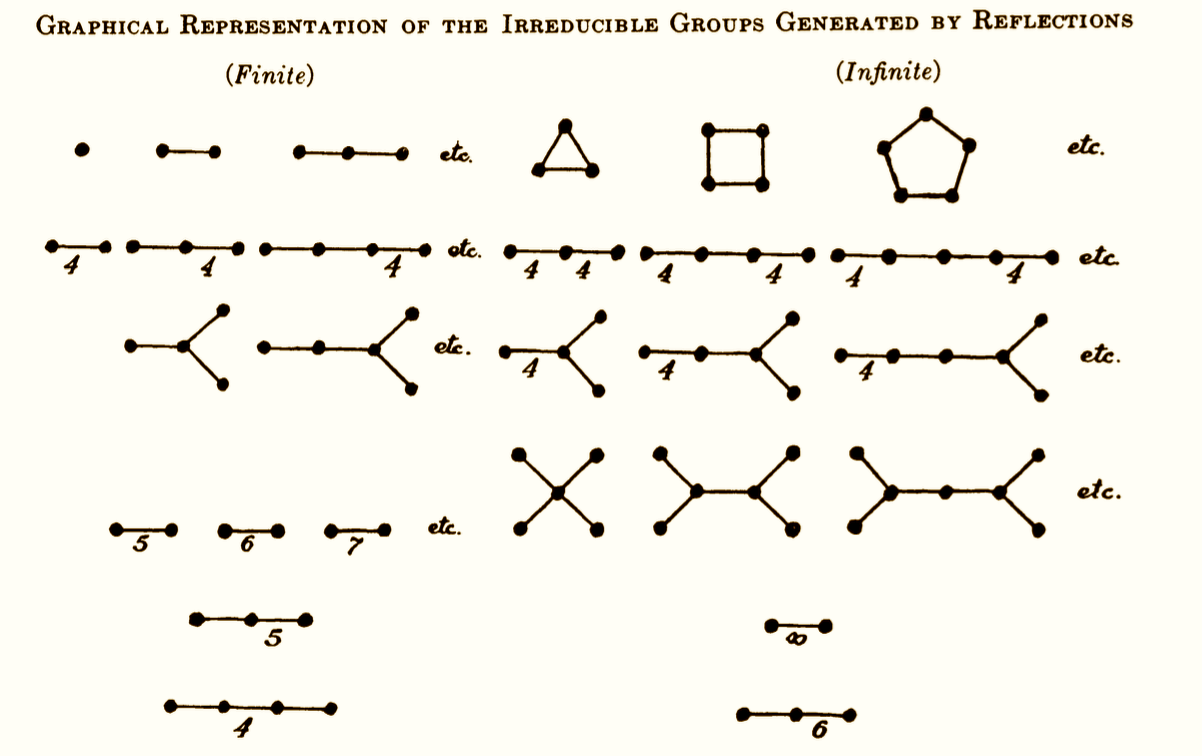} 
 \vspace*{1mm} 
 \caption{Coxeter's pictures of his diagrams} 
 \label{fig:1934-Coxeter}
 \end{figure}

Exactly in the same way in which Du Val introduced in 1934 his dual graphs verbally, without 
drawing them, Coxeter had verbally introduced  in 1931  ``\emph{diagrams of dots and links}'' in order 
to describe the shapes of his spherical simplices (see Figure \ref{fig:1931-Coxeter}). It is only 
in his 1934 paper \cite{C 34} that he published drawings of such graphs 
(see Figure \ref{fig:1934-Coxeter}), which were to be called 
later ``\emph{Coxeter diagrams}'', or ``\emph{Coxeter-Dynkin diagrams}'', 
in reference to their reappearance in a slightly different form in Dynkin's 1946 work \cite{D 46} 
about the structure of Lie groups and Lie algebras.

\begin{figure}%[h!] 
 %\vspace*{6mm}
 \centering 
 \includegraphics[scale=0.55]{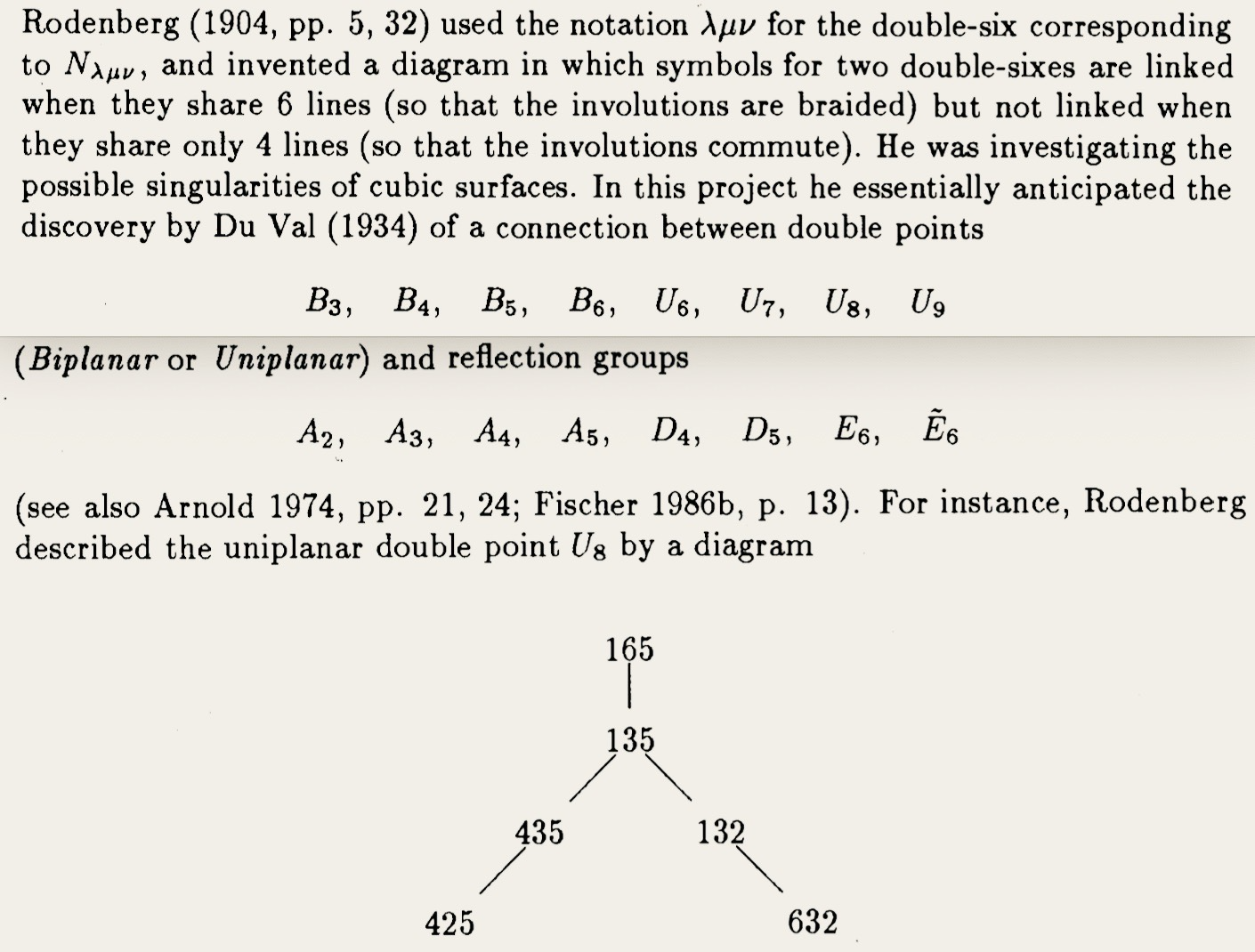} 
 \vspace*{1mm} 
 \caption{One of Rodenberg's diagrams} 
 \label{fig:1904-Rodenberg-3}
 \end{figure}

Much later, Coxeter explained in his 1991 paper \cite{C 91} that analogous diagrams had already 
been introduced by Rodenberg 
in 1904, in his description \cite{R 04} of the plaster models of singular cubic surfaces 
from Schilling's catalog. Rodenberg's intepretation was different, not related to reflections, 
but to special subsets of the configuration of 27 lines on a generic smooth cubic surface 
(see Figure \ref{fig:1904-Rodenberg-3}, containing an 
extract from \cite{C 91}). More details about Rodenberg's convention, based on his 
older paper \cite{R 78}, may be found in Barth and Kn\"orrer's text \cite{BaK 86}.

\section{\bf Mumford's paper on the links of surface singularities}
\label{sec:Mumpap}

One could believe that the combination of Du Val's analogy between his ``trees of curves''  and 
Coxeter's spherical simplices on one side, and Coxeter's diagrams on another side, would 
trigger research on the possible dual graphs of isolated surface singularities. Such an 
interest indeed developed starting from a 1961 paper of David Mumford. 
In this section I explain the aim of Mumford's paper and how it led to the first explicit formulation 
of the notion of dual graph of a configuration of algebraic curves contained in a smooth 
algebraic surface.

\begin{figure}%[h!] 
 %\vspace*{6mm}
 \centering 
 \includegraphics[scale=0.7]{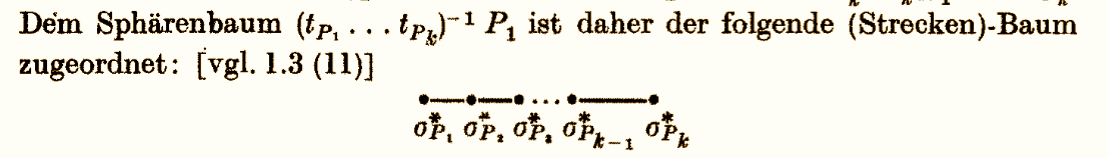} 
 \vspace*{1mm} 
 \caption{The dual graph of Hirzebruch's 1953 resolution paper} 
 \label{fig:1953-Hirzebruch}
 \end{figure}

I could find only one article published between 
1934 and 1961 which contained a drawing of dual graph of 
resolution of isolated surface singularity\footnote{One may think, from a rapid glance, that Du Val's paper \cite{V 44} is an exception. But the graphs of that paper, which I rediscovered with a different 
interpretation in \cite{PP 11}, are not dual graphs of configurations of curves. 
They indicate relations between infinitely 
near points of a given smooth point of an algebraic surface, being variants of 
Enriques' diagrams introduced in \cite{EC 17}.}.  It is Hirzebruch's 
1953 paper \cite{H 53}, in which he proved that one could not only resolve 
the singularities of complex algebraic surfaces, but also of the more 
general \emph{complex analytic} ones. That article contains a single illustration  
(see Figure \ref{fig:1953-Hirzebruch}),  which  
depicts the general shape of possible dual graphs of resolutions 
for a class of singularities which is crucial for his method\footnote{\label{fn:HJsing} Such singularities 
are variously called nowadays ``{\it Hirzebruch-Jung singularities}'', ``{\it cyclic quotient 
singularities}'' or ``{\it toric surface singularities}''. The first name alludes to their 
importance for Hirzebruch's method or resolution of arbitrary complex analytic 
surfaces (explained in \cite{PP 11bis}), inspired by the ideas of Jung's 
method of local parametrization of algebraic surfaces mentioned at the end of 
Section \ref{sec:wmres}. }. 
Unlike the case of Du Val's singularities which do not affect the conditions of adjunction, 
which have resolutions for which all the curves composing 
the exceptional divisor  have self-intersection $-2$, here the self-intersections can  
be arbitrary negative integers\footnote{The ``normal'' surface 
singularities which admit resolutions with such a dual graph, its composing curves 
being moreover smooth, rational and pairwise transversal, are exactly the  
Hirzebruch-Jung singularities alluded to in footnote \ref{fn:HJsing}. 
Their ``\emph{links}'', in the terminology explained later 
in this section, are the so-called \emph{lens spaces}. One may consult Weber's recent paper 
\cite{W 18} for details about the history of the study of lens spaces and their 
relation with Hirzebruch-Jung singularities.}. But a common feature of both cases is that 
all those curves are smooth and rational. This means that from a topological viewpoint they are 
$2$-dimensional spheres. 

Hirzebruch used  
the expression ``\emph{Sph\"arenbaum}'', that is, ``tree of spheres'' 
for the configurations schematically represented in Figure \ref{fig:1953-Hirzebruch}. 
He explained that this terminology had been introduced by Heinz Hopf in his 1951 paper \cite{H 51}  
for the configurations of  $2$-dimensional spheres which are created 
by successive blow-ups of points, starting from a point on a smooth complex algebraic surface. 
Hopf chose that name because those spheres intersect 
in the shape of a tree\footnote{In Section 4 of \cite{H 51}, Hopf wrote: 
``{\it  [...] die $\sigma_i$ lassen sich den 
Eckpunkten eines Baumes -- d. h. eines zusammenh\"angenden Streckencomplexes, der keinen geschlossenen Streckenzug enth\"alt -- so zuordnen, dass zwei $\sigma_i$ dann und nur dann 
einen gemeinsamen Punkt haben, wenn die entsprechenden Eckpunkte eine Strecke begrenzen.}''. 
This is to be compared with 
Hirzebruch's explanation given in section (11) of \cite{H 53}: 
``{\em $K_n$ wird von H. Hopf als 
Sph\"arenbaum bezeichnet, da sich die Sph\"aren $\sigma_i^*$ eineindeutig den Eckpunkten 
eines Baumes zuordnen lassen: Zwei $\sigma_i^*$ haben genau dann einen Schnittpunkt, wenn 
die zugeordneten Eckpunkte im Baum eine Strecke begrenzen.}''}. 

\begin{figure}%[h!] 
 %\vspace*{6mm}
 \centering 
 \includegraphics[scale=0.6]{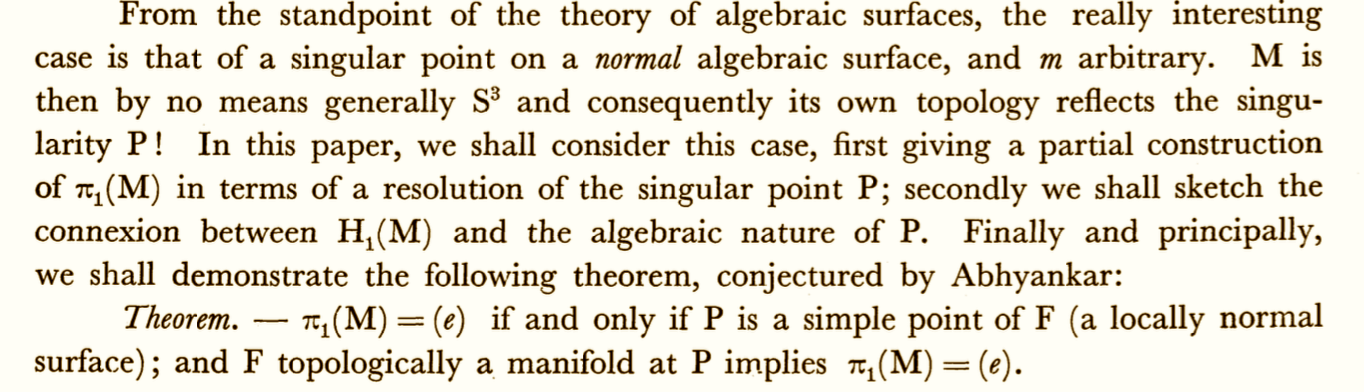} 
 \vspace*{1mm} 
 \caption{Mumford's motivation} 
 \label{fig:1961-Mumford-1}
 \end{figure}
 
 Let us look again at the models of algebraic surfaces from Figures 
 \ref{fig:2018-Hauser}, \ref{fig:2018-Fruhbis-Kruger-1} and \ref{fig:1903-Schilling}. 
 In each case, one has a representation of only part of 
 the surface, as the whole surface is unbounded. The chosen part is obtained 
 by considering the intersection of the entire surface with a ball centered at 
 the singular point under scrutiny. By this procedure, one obtains a 
 portion of the surface possessing a boundary curve. When the ball's radius 
 is small enough, one gets a curve whose qualitative shape (number of connected 
 components, number of singular points on each component, etc.) does not 
 depend on the radius. It is called nowadays the ``\emph{link}'' of the 
 singular point.

 One may perform an analogous construction starting from a point of a \emph{complex} 
 algebraic surface. 
 When the point is taken on a ``\emph{normal}'' complex surface\footnote{Normality is a technical 
 condition which implies that all the singular points of the surface are isolated. One may find 
 details about it in Laufer's book \cite{L 71} or in Mumford's book \cite{M 99}.}, 
 then its associated link is a $3$-dimensional 
 manifold. In the late 1950's, Abhyankar conjectured that it was 
 impossible to obtain a counterexample to the Poincar\'e conjecture following this procedure. 
 In other words, that it was impossible to find a point on a normal complex 
 surface, whose link is simply connected and different from the $3$-dimensional 
 sphere. The aim of Mumford's paper \cite{M 61} was to prove this conjecture, as 
 may be seen in Figure \ref{fig:1961-Mumford-1}, which reproduces the end of its introduction.

\begin{figure}%[h!] 
 %\vspace*{6mm}
 \centering 
 \includegraphics[scale=0.6]{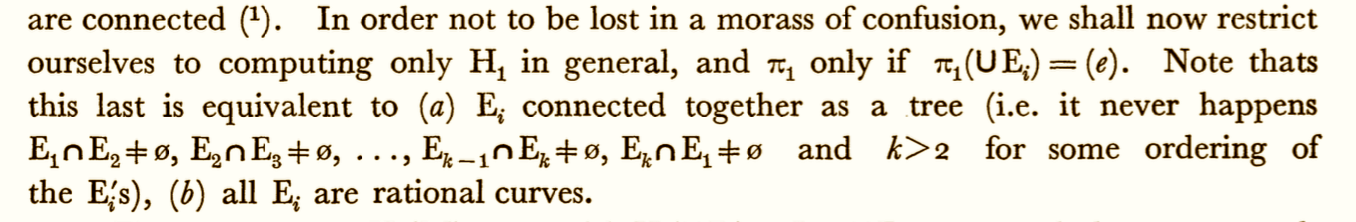} 
 \vspace*{1mm} 
 \caption{Mumford's curves ``connected together as a tree''} 
 \label{fig:1961-Mumford-2}
 \end{figure}

 Mumford's proof started from a resolution of the given point, assuming -- which is 
 always possible -- that its  exceptional divisor has 
 ``\emph{normal crossings}''\footnote{\label{fn:normcross} This 
 means that the exceptional divisor is  composed 
 of smooth curves intersecting pairwise transversally,  and that three such curves do not 
 have common points.}. It proceeded along the following steps\footnote{One needs to 
 know a certain amount of algebraic topology in order to understand this proof completely. 
 The reader with a strong taste for visualization may learn the needed notions from the website 
 \cite{AS 17}.}:
    \begin{enumerate}
        \item \label{step1} Show that the link $M$ of the given point $p$ of a normal complex surface 
           is determined by the exceptional divisor and 
            by the self-intersection numbers of its composing curves.
        \item  Show that if $M$ is simply connected, then the curves of the exceptional divisor 
           are ``\emph{connected together as a tree}'' and are all rational 
           (see Figure \ref{fig:1961-Mumford-2}). 
        \item \label{step3} Under this hypothesis, write a presentation of the fundamental group $\pi_1(M)$ 
           of $M$ in terms of the configuration of curves of the exceptional divisor and of 
           their self-intersection numbers. 
        \item Deduce from this presentation that if $\pi_1(M)$  is trivial, then 
           one may contract algebraically 
           one of the curves of the exceptional divisor to a point and obtain again  a 
           resolution whose exceptional divisor has normal crossings. 
        \item Iterating such contractions, show that the given point $p$ is a smooth point 
            of the starting surface, which implies that its link $M$ is the $3$-dimensional sphere. 
    \end{enumerate}
    
    \begin{figure}%[h!] 
 %\vspace*{6mm}
 \centering 
 \includegraphics[scale=0.6]{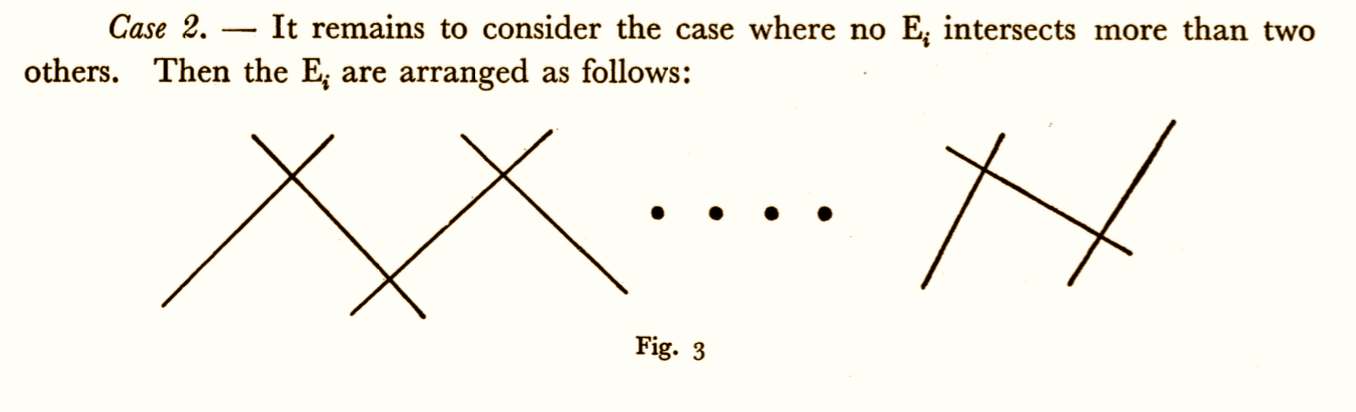} 
 \vspace*{1mm} 
 \caption{Mumford's diagram} 
 \label{fig:1961-Mumford-diagr}
 \end{figure}

In what concerns step (\ref{step1}), Mumford proved in fact 
that the link $M$ is determined by the dual graph 
of the exceptional divisor, decorated by the genera and the self-intersection 
numbers of the associated curves. This formulation is slightly anachronistic,  
because he still did not formally introduce this dual graph. He said only that the curves 
of the exceptional divisor were ``\emph{connected together as a tree}'', 
which is similar to Du Val's terminology of his 1934 paper discussed in Section \ref{sec:DVCox} 
(see again Figure \ref{fig:1934-Du-Val-graphes-duaux}). Unlike 
Hirzebruch in his 1953 article, he did not even draw a dual graph, 
but only a schematic representation of the same type of configuration of curves as that 
of Hirzebruch's paper \cite{H 53} (see Figure \ref{fig:1961-Mumford-diagr}). One may 
notice that the same drawing convention was to be followed by Michael Artin in his 1962 paper \cite{A 62} 
(see the upper half of Figure \ref{fig:1962-6-Artin}), before his switch to dual graphs in 
the 1966 paper \cite{A 66}.

\begin{figure}%[h!] 
 %\vspace*{6mm}
 \centering 
 \includegraphics[scale=0.5]{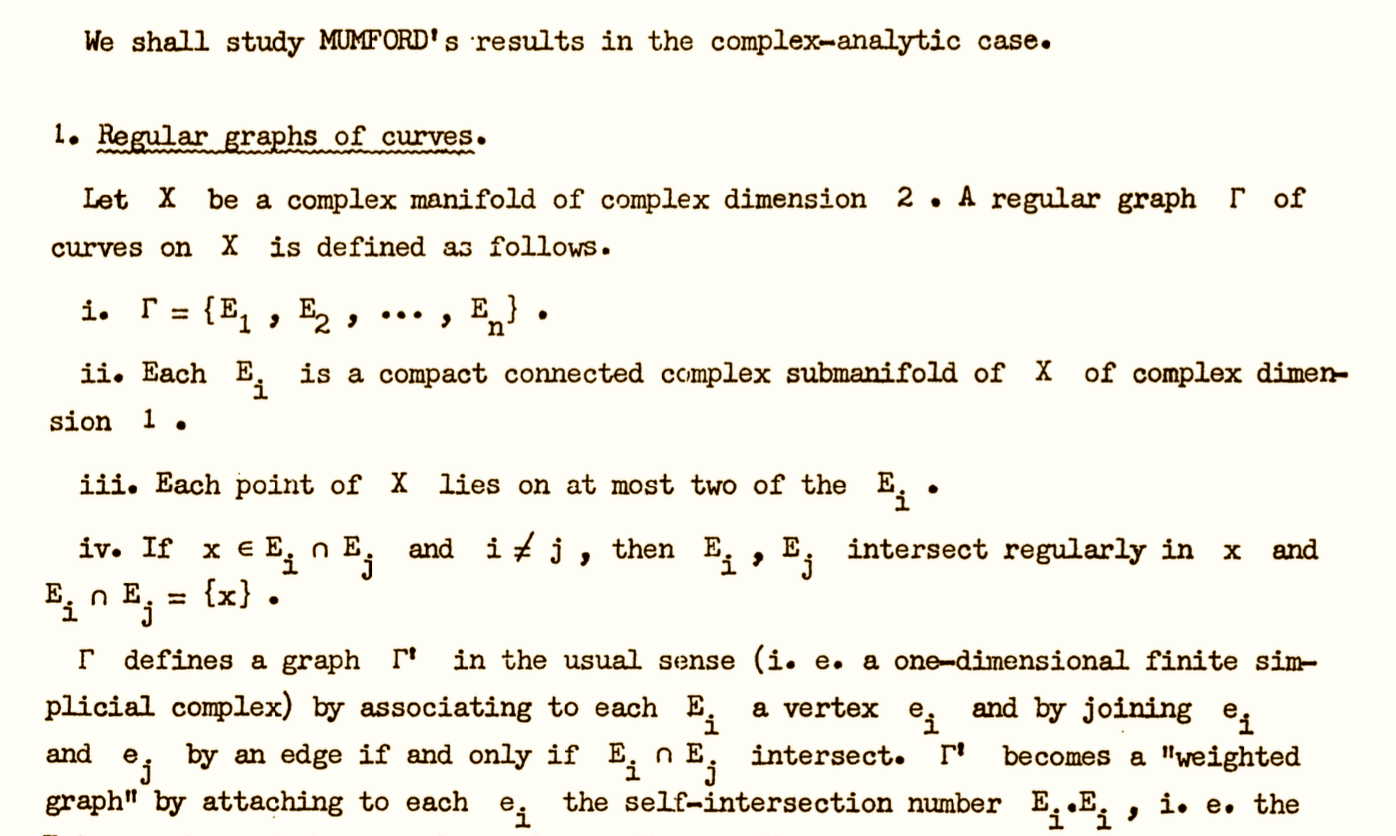} 
 \vspace*{1mm} 
 \caption{Hirzebruch's switch from ``graphs of curves'' to graphs ``in the usual sense''} 
 \label{fig:1962-Hirzebruch}
 \end{figure}

 It seems that the notion of dual graph of an arbitrary curve configuration, not necessarily formed 
 of smooth rational curves, was formally defined for the first time
 in Hirzebruch's 1963 Bourbaki Seminar talk \cite{H 63} discussing the previous results of Mumford.  
 Hirzebruch associated to any ``regular graph of curves'' a ``graph in the usual sense'' 
 (see Figure \ref{fig:1962-Hirzebruch}), without using 
 the terminology ``dual graph'', which seems to have appeared later in the 1960s. 
 Then, he stated the result of step (\ref{step3}) formulated above as the fact that the 
 dual graph determines an explicit presentation of the fundamental group of the link $M$ 
 -- of course, under Mumford's hypothesis that the exceptional divisor has normal 
 crossings, that all its components are rational curves and that the graph is a tree. 
 
 Less than 10 years later, appeared the first books explaining -- among other things -- algorithms 
 for the computation of dual graphs of resolutions of singularities of normal surfaces: 
 Hirzebruch, Neumann and Koh's book \cite{HNK 71} and Laufer's book \cite{L 71}. 
 All those algorithms followed Hirzebruch's method of his 1953 paper from which  
 Figure \ref{fig:1953-Hirzebruch} was extracted.

Before passing to the next section, let me quote an e-mail received on 9 January 2018, 
in which Mumford answered my questions about the evolution of the notion of dual graph: 

\begin{quote}
     ``Perhaps the following is useful. In much of the 20th century, math papers never 
     had any figures. As a geometer, I always found this absurd and frustrating. 
     In my ``red book'' intro to AG [Algebraic Geometry, \cite{M 99}], 
     I drew suggestive pictures of various schemes, 
     trying to break through this prejudice. On the other hand, I listened to many lectures 
     by Oscar Zariski and, on rare occasions, we, his students, noticed him making a small 
     drawing on the corner of the blackboard. You see, the Italian school had always in mind 
     actual pictures of the real points on varieties. Pictures of real plane curves and plaster casts 
     of surfaces given by the real points were widespread. If you want to go for firsts, check out 
     Isaac Newton's paper classifying plane cubics.  So we were trained to ``see'' the resolution 
     as a set of curves meeting in various ways. The old Italian theory of ``infinitely near points'' 
     was, I think, always drawn that way. Of course, this worked out well for compactifying 
     moduli space with stable curves.  I'm not clear who first changed this to the dual graphs. 
     Maybe it was Fritz [Hirzebruch].'' 
\end{quote}

\section{\bf Waldhausen's graph manifolds and Neumann's calculus with graphs}
\label{sec:WN}

Mumford's theorem stating that the link of a complex normal surface singularity is determined by 
the dual graph of any of its resolutions whose exceptional divisor has normal crossings 
 raised the question whether, conversely, it was possible to recover the 
 dual graph from the structure of the link. 
 
 Formulated in this way, the problem cannot be solved, because resolutions are not unique. 
 Indeed, given a resolution whose exceptional divisor has normal crossings, one can get 
 another one by blowing up any point of the exceptional divisor. The new resolution has a 
 different dual graph, with one additional vertex. Is there perhaps a 
 minimal resolution, from which all other resolutions are obtained by sequences of blow ups 
 of points? Such a resolution indeed exists\footnote{The figures  in 
 Sections \ref{sec:intro} and \ref{sec:meangraph} present in fact dual 
 graphs of \emph{minimal} resolutions of the corresponding surface singularities.}, 
 and one may ask instead whether its dual graph is 
 determined by the corresponding link.

This second question was answered affirmatively in the 1981 paper \cite{N 81} of Walter Neumann, 
building on a 1967 paper \cite{W 67} of Friedhelm Waldhausen. Let me describe successively 
the two papers, after a supplementary discussion of Mumford's article \cite{M 61}. 

Mumford looked at the link $M$ of a singular point of a normal complex algebraic surface 
as the boundary of a suitable ``tubular neighborhood'' of the exceptional divisor of the chosen 
resolution. One of his crucial insights was that the assumption that this divisor 
has normal crossings\footnote{Recall that this 
notion was explained in footnote \ref{fn:normcross}.} implies that the $3$-dimensional manifold $M$  
may be described only in terms of real curves and surfaces, which are objects of 
smaller dimension. Namely, the link $M$ may be obtained by suitably cutting and pasting  
continuous families of circles -- called {\it circle bundles} -- 
parametrized by the points of the curves of the exceptional divisor. Let me explain why.

In the simplest case where 
the exceptional divisor is a single smooth algebraic curve -- topologically a {\it real surface}, because 
it is a \emph{complex  curve} -- 
such a tubular neighborhood has a structure of disk bundle over the curve: one may fill it by 
discs transversal to the curve. 
Being its boundary, the link has therefore a structure of circle bundle over this surface. 

In general, the exceptional divisor has several components. Then the boundary of one of its  
tubular neighborhoods has again a structure 
of circle bundle far from the intersection points of those components, but one has to make 
a careful analysis near such points. Mumford worked in special  
neighborhoods of them, which he called ``\emph{plumbing fixtures}". 
These allow to see in which way one passes 
from the circle bundle over a curve of the exceptional divisor to that over a second such  
curve, intersecting the first one at the chosen point. 
In terms of the dual graph, the description is very simple: 
to every edge of it one can assign a ``plumbing fixture''. In the link, which is identified with 
the boundary of the tubular neighborhood, it gives rise to a torus. 
By moving inside the link $M$ and crossing this torus, 
 one passes from the first circle fibration to the second one. In order 
to understand precisely in which way the transition is made, one has to look at the relative 
positions of the circles of both fibrations on the separating torus. 
These intersect transversally at exactly one point. 

\begin{figure}%[h!] 
 %\vspace*{6mm}
 \centering 
 \includegraphics[scale=0.6]{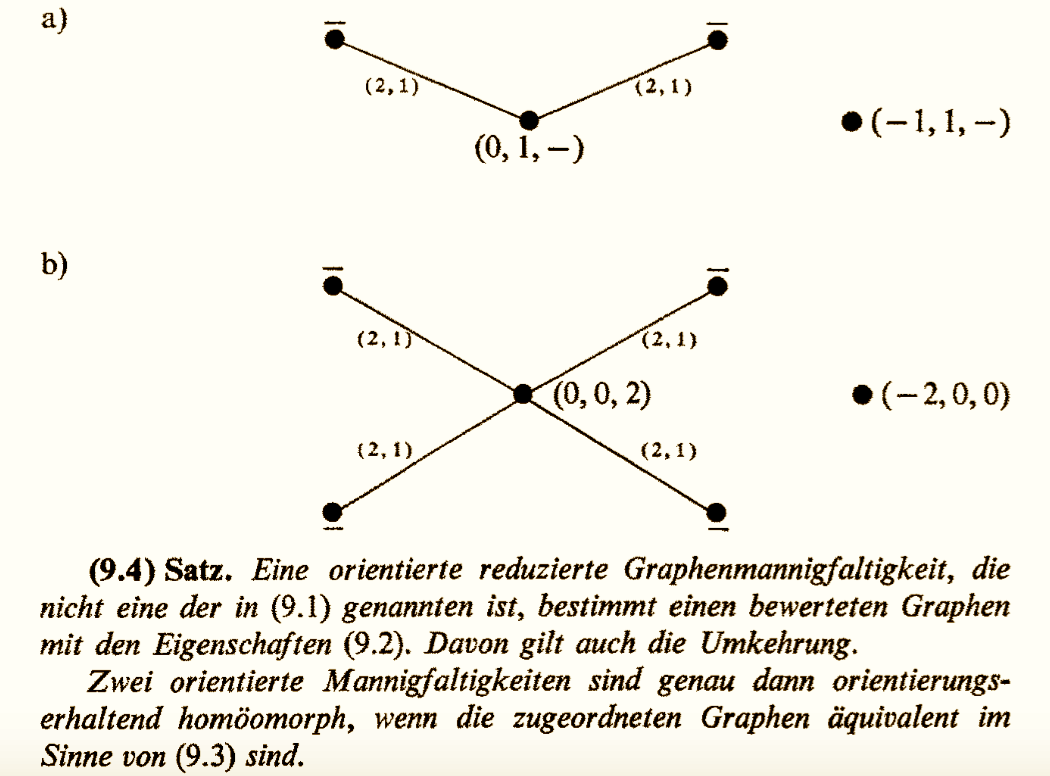} 
 \vspace*{1mm} 
 \caption{Two of Waldhausen's ``graph manifolds''} 
 \label{fig:1967-Waldhausen}
 \end{figure}

This phenomenon gave rise to the notion of ``\emph{plumbed $3$-manifold}''. It is a 
$3$-dimensional manifold 
constructed from a decorated graph by performing a ``\emph{plumbing}'' operation for each 
edge of the graph, similar to that described by Mumford in his paper. 
Each vertex comes equipped with two numbers, one 
representing the genus of a surface and the second its self-intersection number in an associated 
disk-bundle of dimension four.  There is a subtlety related 
to orientations, which obliges one to decorate the edges with signs. 

One witnesses here a metamorphosis of the interpretation of the weighted dual graphs. 
If they started by representing the configurations of curves obtained as exceptional 
divisors of resolutions of singularities, they became blueprints for building  
certain $3$-manifolds. It was Waldhausen who developed a subtle theory 
of those manifolds in \cite{W 67}. He called them ``\emph{Graphenmannigfaltigkeiten}'' 
-- that is, ``\emph{graph-manifolds}'' 
and not ``plumbed manifolds'', in order to emphasize the idea that they are defined by graphs. 
In fact, he considered slightly more general graphs, whose edges are also decorated with 
pairs of numbers (see Figure \ref{fig:1967-Waldhausen}). This convention allowed the transitions 
from one circle fibration to another one across a torus to be performed by letting the 
fibers from both sides intersect in any way, not necessarily transversally at a single point. 
One of his main theorems states that any graph-manifold has a unique minimal 
graph-presentation, except for an explicit list of ambiguous manifolds. 

\begin{figure}%[h!] 
 %\vspace*{6mm}
 \centering 
 \includegraphics[scale=0.5]{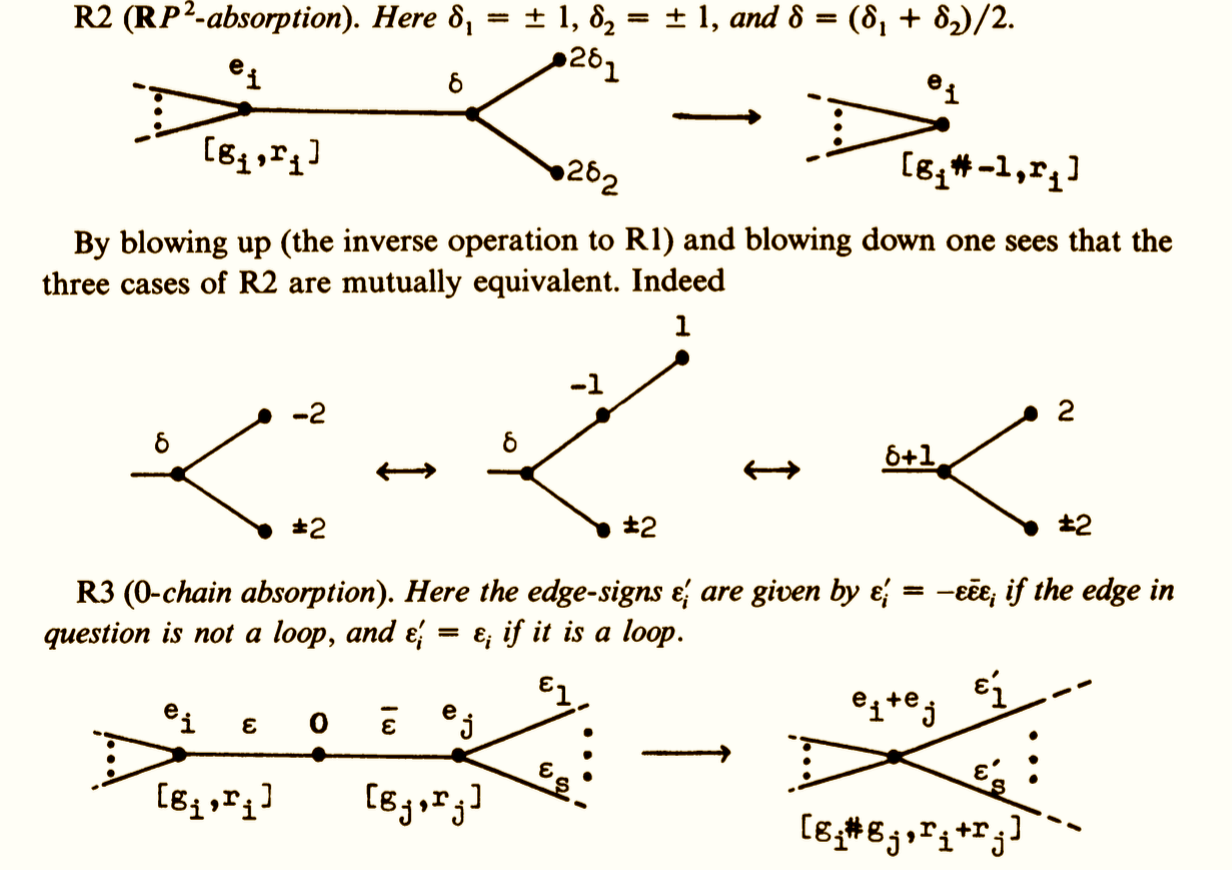} 
 \vspace*{1mm} 
 \caption{Part of Neumann's ``plumbing calculus''} 
 \label{fig:1981-Neumann}
 \end{figure}

In his 1981 paper \cite{N 81}, Neumann turned this theorem into an  
algorithm, allowing to determine whether two weighted graphs in the original sense of the 
plumbing operations determined the same oriented $3$-dimensional manifold. 
Roughly speaking, this algorithm 
consists in applying successively the rules of a ``\emph{plumbing calculus}'' -- some of them being 
represented in Figure \ref{fig:1981-Neumann} -- in the direction which diminishes the 
number of vertices of the graph. Two graphs determine the same $3$-dimensional 
manifold if and only if the associated ``minimal'' graphs coincide -- again, up to a little 
ambiguity related to the signs on the edges.  

The algorithm allowed Neumann to prove important 
topological properties of normal surface singularities and of families of smooth complex 
curves degenerating  to singular ones. For instance, he showed that the decorated dual 
graph of the minimal resolution with normal crossings is determined by the 
\emph{oriented} link of the singularity. This unified the two viewpoints 
on the graphs associated to surface singularities discussed in this paper (as dual graphs of their resolutions, and as blueprints for building their links).

\section{\bf  Conclusion}
\label{sec:concl}

We could continue this presentation of the interaction between graphs and singularities  
in several directions: 

\begin{itemize}
    \item By examining other classification problems of singularity theory which led to 
        lists of dual graphs, between Mumford's paper \cite{M 61} and Neumann's paper 
        \cite{N 81} (for instance Brieskorn's paper \cite{B 68}, 
        Wagreich's paper \cite{W 72} and Laufer's papers \cite{L 73, L 77}). 
        We saw that it was such a classification problem which led Du Val to 
        his consideration of ``graphs of curves'', and which led to Artin's lists 
        of dual graphs shown in Figures \ref{fig:1962-6-Artin}  and \ref{fig:1966-Artin-rtp}. 
        Another such problem led to the more recent paper of Chung, Xu and Yau 
        from which Figure \ref{fig:2009-Chung-Xu-Yau} is extracted. Note that in \cite{L 73}, 
        Laufer classified all ``\emph{taut}'' normal surface singularities, that is, 
        those which are determined up to 
        complex analytic isomorphisms by the dual graphs of their minimal resolutions.  
        He showed in particular that all rational double and triple points of Figures 
        \ref{fig:1962-6-Artin} and \ref{fig:1966-Artin-rtp} are taut. 
        This result had already been proved by Brieskorn in \cite{B 68} for rational double points.  
        As a consequence, this class of singularities coincides with Du Val's 
        singularities which do not affect the condition of adjunction. This result is much stronger 
        than the fact that their minimal resolutions have the same dual graphs.
   \item By examining the applications and developments of Neumann's ``plumbing 
      calculus''. One could analyze its variant developed by Eisenbud and 
      Neumann in \cite{EN 85} for the study of certain links (that is, disjoint unions 
      of knots) in integral homology 
      spheres which are graph-manifolds, its applications initiated by Neumann \cite{N 89}  
      to the study of complex plane curves at infinity, 
      or those initiated by N\'emethi and Szilard \cite{NS 12} and continued by Curmi 
      \cite{C 20} to the study of boundaries 
      of ``\emph{Milnor fibers}'' of non-isolated surface singularities. 
   \item By discussing generalizations of dual graphs to higher dimensions. In general, 
      when one has a configuration of algebraic varieties, one may represent them by points, 
      and fill any subset of the total set of such points by a simplex, whenever the corresponding 
      varieties have a non-empty intersection. One gets in this way the so-called 
      ``\emph{dual complex}'' of the configuration of varieties. In the same way as there was a 
      substantial lapse of time since the idea of dual graph emerged till it became an active object 
      of study, an analogous phenomenon occurred with this more general notion. 
      It seems to have appeared independently in the 1970s, in Danilov's paper \cite{D 75} 
      -- whose results were rediscovered with a completely different proof 
      by Stepanov in his 2006 article \cite{S 06} -- in Kulikov's paper \cite{K 77}  
      and in Persson's book \cite{P 77}. 
      Information about recent works on dual complexes may be found in 
      Payne's paper \cite{P 13}, in Koll\'ar's paper 
      \cite{K 13} and in the paper \cite{FKX 17} by de Fernex, Koll\'ar and Xu.  
      One may use Nicaise's paper \cite{N 16} as an introduction to the relations  between 
      dual complexes and 
      ``\emph{non-Archimedean analytifications in the sense of Berkovich}''. 
   \item By presenting the notion of ``\emph{fan}'' of the divisor at infinity of a toroidal variety, 
       introduced by Kempf, Knudson, Mumford and Saint-Donat in the 1973 book \cite{KKMS 73}. 
       It is a complex of cones associated to special kinds of configurations of hypersurfaces 
       in complex algebraic varieties. When the configuration has normal crossings, 
       the projectivisation of the fan is in fact the dual complex of the configuration.  
       Fans had been introduced before by Demazure for ``\emph{toric varieties}'' in the 1970 paper 
       \cite{D 70}, and since then they were mainly used in ``\emph{toric geometry}''. 
       Following this direction, we could 
       arrive at the notion of ``\emph{geometric tropicalization}'', which expresses 
       ``\emph{tropicalizations}'' of 
       subvarieties of algebraic tori in terms of the dual complexes of the divisors at infinity 
       of convenient compactifications (see \cite{HKT 09} and \cite[Theorem 6.5.15]{MS 15}). Note that 
       Berkovich's analytification (alluded to at the end of the previous item) and 
       tropicalization are intimately related, as explained by Payne in \cite{P 09, P 15}. 
       Note also that the paper \cite{GBGPPP 18} studies dual graphs 
       of resolutions of normal surface singularities in the same spirit. 
   \item By discussing how Waldhausen's theory of graph-manifolds led to Jaco-Shalen-Johannson's 
      theory of canonical decompositions of arbitrary orientable and closed $3$-manifolds into elementary 
      pieces, by cutting them along spheres and tori (see Jaco and Shalen's book \cite{JS 79} 
      and Johansson's book \cite{J 79}). This is turn gave rise to Thurston's 
      geometrization conjecture of \cite{T 82}, proved partially by Thurston, and which 
      was finally completely settled by Perelman's work \cite{P 03}. For details on Perelman's 
      strategy, one may consult the monographs \cite{BBBM 10} of Bessi\`eres, Besson, Boileau, 
       Maillot and Porti and \cite{MT 14} of Morgan and Tian. 
   \item By speaking about the second, more recent, main source of graphs in singularity theory:  
        the dual graphs of configurations of ``\emph{vanishing cycles}'' in Milnor fibers of isolated 
        hypersurface singularities. Such dual graphs, called sometimes ``\emph{Dynkin diagrams}'', 
        began to be described and drawn after 1970 for special classes of singularities by 
        A'Campo \cite{AC 75, AC 75bis}, Gabrielov \cite{Ga 73, Ga 74} 
        and Gusein-Zade \cite{G 74}. One may consult Arnold's papers \cite{A 73, A 75},  
        Gusein-Zade's survey \cite{G 77} and Brieskorn's papers \cite{B 81, B 83} 
        for a description of the context 
        leading to those researches on Dynkin diagrams and of their relations with other 
        invariants of hypersurface singularities. Du Val's singularities 
        possess configurations of vanishing cycles isomorphic to the dual graphs of their 
        minimal resolutions indicated in Figure \ref{fig:1962-6-Artin}. One may consult 
        Brieskorn's paper \cite{B 00} for a description of the way he proved this theorem instigated  
        by a question of Hirzebruch. In fact, this property characterizes Du Val's singularities (see Durfee's 
        survey \cite{D 79} of many other characterizations of those singularities). In general, 
        the relation between the two types of dual graphs is still mysterious. 
        Note that Arnold described in \cite{A 75} a ``\emph{strange duality}'' inside 
        a set of 14 ``\emph{exceptional unimodular singularities}'', relating the two types of dual graphs. This 
        duality was explained by Pinkham \cite{Pi 77} on one side and 
        Dolgachev and Nikulin \cite{DN 77} on another side (see Dolgachev's Bourbaki seminar 
        presentation \cite{D 83}). Later, Dolgachev related 
        it in \cite{D 96} to the very recent phenomenon -- at the time -- of ``\emph{mirror symmetry}'', 
        but this seems to be only the tip of an iceberg. 
  \end{itemize}

 I will not proceed in such directions, because this would be very difficult to do while remaining 
 reasonably non-technical. I made nevertheless the previous list in order to show that 
 dual graphs and their generalizations to higher dimensions are nowadays common 
 tools in singularity theory, in algebraic geometry and in geometric topology. 
 It is for this reason that I found interesting to examine their births and their early uses. 
 
 We saw that dual graphs of surface singularities were first used mainly verbally, in expressions 
 like ``tree of curves'', ``Sph\"arenbaum''. Drawing them became important for stating results 
 of various problems of classification. This made their verbal description first 
 too cumbersome, then completely inadequate for the description of the wealth of 
 morphologies under scrutiny. 
 Then, their reinterpretation as blueprints for building graph-manifolds led to the 
 development of a ``plumbing calculus'', which transformed them into objects of algebra. 
 The necessity to develop an analogous ``calculus'' appears every time one gets 
 many different encodings of the structure of an object, leading to the problem 
 of deciding which encodings correspond to the same object (another instance of this 
 phenomenon is Kirby's calculus of \cite{K 78}). In other situations 
 -- for instance, that of finite presentations of discrete groups -- it is known that 
 the problem is undecidable. But for plumbing graphs it is solvable, as shown 
 by the works of Waldhausen and Neumann alluded to before.

\medskip
\end{document}